\documentclass[10pt]{article}
  \usepackage{amsfonts, amssymb, amsmath}
  \usepackage{calrsfs, pifont}
  \usepackage{makeidx}
  \usepackage{latexsym}
  \usepackage{mathrsfs}
  \usepackage{bbold}
  \usepackage{indentfirst}
 \def\diag{\mathop{\fam0 diag}}
 \def\End{\mathop{\fam0 End}}
 \def\On{\mathop{\fam0 On}}
 \def\upwardarrow{\mathord{
  \hbox to 5pt{\hss$\vcenter{\hbox to 2.4pt{\hss$\mathchar"222$\hss}\hrule}\hss$}
 }}
 \def\downwardarrow{\mathord{
  \hbox to 5pt{\hss$\vcenter{\hrule\hbox to 2.4pt{\hss$\mathchar"223$\hss}}\hss$}
 }}
 \def\to{\rightarrow}
 \def\osum{o\text{-}\!\sum}
 \def\oleq{\mathop{{\bigcirc\hspace{-8.5pt}\raisebox{1pt}{$\scriptstyle\leq$}\hspace{2pt}}}}
 \def\beginproof{\smallskip\par\text{$\vartriangleleft$}}
 \def\endproof{\text{$\vartriangleright$}}
 \newcommand{\subsubsec}[1]{\qquad{\bf #1}}
 \newcommand{\im}{\mathop{\fam0 Im}}
 \def\Orth{\mathop{\fam0 Orth}\nolimits}
 \def\Hom{\mathop{\fam0 Hom}\nolimits}
 \def\mix{\mathop{\fam0 mix}}

 \newcommand{\dom}{\mathop{\fam0 dom}}

 \newcommand{\Theorem}[1]{\smallskip\textbf{{Theorem#1}.\/ }\sl}
 \newcommand{\theorem}[1]{{\textbf{Theorem#1.\/~}}\sl}
 \newcommand{\Corollary}[1]{\smallskip\textbf{{Corollary#1}.\/~}\sl}
 \newcommand{\corollary}[1]{\textbf{{Corollary#1}.\/~}\sl}

 \newcommand{\proclaim}[1]{{\textbf{#1\/~}}\sl}
 \newcommand{\Endproc}{\rm}

\begin{document}
 \begin{center}
 \large \textbf{BOOLEAN VALUED ANALYSIS \\OF ORDER BOUNDED OPERATORS}
 \end{center}
 \begin{center}
 A.~G.~Kusraev and S.~S. Kutateladze
 \end{center}
 \begin{abstract}
 This is a~survey of some recent applications of  Boolean valued models of set theory
 to  order bounded operators in vector lattices.\\[6pt]
 \end{abstract}
\markboth
{\scshape Boolean Valued Analysis of Order Bounded Operators}
{\scshape A.~G.~Kusraev and S.~S. Kutateladze}

 \section{Introduction}

 The term {\it Boolean valued analysis} signifies the technique of studying properties of an
 arbitrary mathematical object by comparison between its representations in two different
 set-theoretic models whose construction utilizes principally distinct Boolean algebras.
 As these models, we usually take the classical Cantorian paradise in the shape of the von
 Neumann universe and a specially-trimmed Boolean valued universe in which the conventional
 set-theoretic concepts and propositions acquire  bizarre interpretations. Use of two models
 for studying a single object is a family feature of the so-called {\it nonstandard methods of analysis}.
 For this reason, Boolean valued analysis  means an instance of nonstandard analysis in common parlance.

 Proliferation of Boolean valued analysis stems from the celebrated achievement of P.~J.~Cohen who proved
 in the beginning of the 1960s that the negation of the continuum hypothesis, $\mathrm{CH}$, is consistent
  with the  axioms of Zermelo--Fraenkel set theory, $\mathrm{ZFC}$. This result by~Cohen, together
  with the consistency of $\mathrm{CH}$ with $\mathrm{ZFC}$ established earlier  by K.~G\"odel, proves
  that $\mathrm{CH}$ is independent of the conventional axioms of $\mathrm{ZFC}$.

 The first application of Boolean valued models to functional analysis
 were given by E.~I.~Gordon for Dedekind complete vector lattices and positive operators in
 \cite{Gor1}--\cite{Gor3} and  G.~Takeuti self-adjoint operators in Hilbert spaces and harmonic analysis
 in \cite{Tak}--\cite{Tak1}. The further
 developments and corresponding references are presented in \cite{BVA, IBA}.

 The aim of the paper is to survey some recent applications of  Boolean valued models of set theory
 to studying order bounded operators in vector lattices. Chapter~1 contains a sketch of the adaptation to analysis of
  the main
 constructions and principles of Boolean valued models of set theory. The three
 subsequent chapters treat the classes of operators in vector lattices:
  multiplication type operators, weighted shift type operators, and conditional expectation type operators.

 The~reader can find the~necessary information on Boolean algebras in \cite{Sik, Vl}, on the~theory
 of vector lattices, in~\cite{AB, KVP, DOP, Vul, Z}, on
 Boolean valued models of set theory, in \cite{Bell, Jech, TZ}, and on
 Boolean valued analysis, in~\cite{KK, BVA, IBA}.

 Everywhere below $\mathbb{B}$ denotes a complete Boolean algebra, while $\mathbb{V}^{(\mathbb{B})}$
 stands for the corresponding Boolean valued universe (the universe of $\mathbb{B}$-valued sets).
 A~{\it partition of~unity\/} in~$\mathbb{B}$ is a~family
 $(b_{\xi})_{\xi\in \Xi}\subset\mathbb{B}$ with
 $\bigvee_{\xi\in \Xi} b_{\xi}=\mathbb{1}$ and
 $b_{\xi}\wedge b_{\eta}=\mathbb{0}$ for $\xi\ne \eta$.

 By a~vector lattice  throughout the~sequel we will mean
 a~real Archimedean vector lattice, unless specified otherwise. We~let $\!:=$ denote
 the~assignment by definition, while
 $\mathbb{N}$, $\mathbb{Q}$, $\mathbb{R}$, and $\mathbb{C}$ symbolize the~naturals,
 the~rationals, the~reals, and the~complexes. We denote  the~Boolean algebras of bands
 and band projections in a~vector lattice~$X$  by  $\mathbb{B}(X)$ and $\mathbb{P}(X)$;
 and we let $X^\mathrm{u}$ stand for
 the~universal completion of a vector lattice $X$.

 The ideal center $\mathcal{Z}(X)$ of a vector lattice $X$ is an $f$-algebra. Let $\Orth(X)$ and $\Orth^\infty(X)$ stand for the $f$-algebras of orthomorphisms and extended orthomorphisms,
 respectively $X$. Then $\mathcal{Z}(X)\subset\Orth(X)\subset\Orth^\infty(X)$. The space of all order bounded linear operators from $X$ to $Y$ is denoted by $L^\sim(X,Y)$. The Riesz--Kantorovich Theorem tells us that if $Y$ is a Dedekind complete vector lattice then so is $L^\sim(X,Y)$.

 \section{Chapter 1. Boolean Valued Analysis}

 \subsection{1.1.~Boolean Valued Models}

 We start with recalling some auxiliary facts about the
 construction and  treatment of Boolean valued models.
 Some more detailed presentation can be found in~\cite{Bell, BVA, IBA}. In the sequel
 $\mathrm{ZFC}\!:=\mathrm{ZF}+\mathrm{AC}$, where $\mathrm{ZF}$ stands for
 the Zermelo--Fraenkel set theory and $\mathrm{AC}$ for the axiom of
 choice.

 \subsubsec{1.1.1.}~Let ${\mathbb B}$ be a~complete Boolean algebra. Given an ordinal
 $\alpha$, put
 $$
 {\mathbb V}_{\alpha}^{({\mathbb B})}\!:=\bigl\{\,x\,:\,x \mbox{
 is\ a\ function}\ \wedge\ (\exists\,\beta)\bigl(\beta<\alpha\ \wedge\
 \dom (x)
 \subset {\mathbb V}_{\beta}^{({\mathbb B})}\ \wedge\ \im
 (x)\subset {\mathbb B}\bigr)\bigr\}.
 $$

 After this recursive definition the {\it Boolean valued
 universe\/} ${\mathbb V}^{({\mathbb B})}$ or, in other words, the
 {\it class of ${\mathbb B}$-sets\/} is introduced by
 $$
 {\mathbb V}^{({\mathbb B})}\!:=\bigcup\limits_{\alpha\in\On}
 {\mathbb V}_{\alpha}^{({\mathbb B})},
 $$
 with $\On$ standing for the class of all ordinals.

 In case of the two-element Boolean algebra
 $\mathbb 2\!:= \{\mathbb{0}, \mathbb 1\}$ this procedure yields a version of the classical
 {\it von Neumann universe\/} ${\mathbb V}$ where
 $\mathbb{V}_0\!:=\varnothing$,
 $\mathbb{V}_{\alpha+1}\!:=\mathcal{P}(\mathbb{V}_{\alpha})$,
 $\mathbb{V}_{\beta}\!:=\bigcup\nolimits_{\alpha<\beta}\mathbb{V}_{\alpha}$,
 $\beta$ is a limit ordinal (cp.~\cite[Theorem 4.2.8]{IBA}).

 \subsubsec{1.1.2.}~Let $\varphi$ be an arbitrary formula of $\mathrm{ZFC}$, Zermelo--Fraenkel set theory with  choice. The {\it Boolean truth value\/} $[\![\varphi]\!]\in {\mathbb B}$ is introduced by induction on the complexity of~$\varphi$ by naturally interpreting the propositional connectives and quantifiers in the Boolean algebra ${\mathbb B}$
 \big(for instance, $[\![\varphi_1\lor\varphi_2]\!]\!:=[\![\varphi_1]\!]\lor[\![\varphi_2]\!]$\big)
 and taking into consideration the way in which a formula is built up from atomic formulas. The Boolean truth values of the {\it atomic formulas\/} $x\in y$ and $x=y$ \big(with $x,y$ assumed to be the elements of ${\mathbb V}^{({\mathbb B})}$\big) are defined by means of the following recursion schema:
 $$
 \gathered {} [\![x\in y]\!]= \bigvee\limits_{t\in\dom(y)}\!\!
 \bigl(y(t)\wedge [\![t=x]\!]\bigr),
 \\%
 [0.5\jot]
 [\![x=y]\!]=\bigvee\limits_{t\in\dom(x)}\!\!\bigl(x(t)\Rightarrow [\![t\in
 y]\!]\bigr)\wedge \bigvee\limits_{t\in\dom(y)}\!\!\bigl(y(t)\Rightarrow [\![t\in
 x]\!]\bigr).
 \endgathered
 $$
 The sign $\Rightarrow$ symbolizes the implication in ${\mathbb
 B}$; i.\,e., $(a\Rightarrow b)\!:=(a^\ast\vee b)$, where $a^\ast$ is as
 usual the {\it complement\/} of~$a$.
 The universe ${\mathbb V}^{({\mathbb B})}$ with the  Boolean
 truth value of a~formula is a~model of set theory in the sense
 that the following  is fulfilled:

 \subsubsec{1.1.3.} \proclaim{Transfer Principle.}For every theorem
 $\varphi$ of~$\mathrm{ZFC}$, we have
 $[\![\varphi]\!]=\mathbb 1$ {\rm(}also in~$\mathrm{ZFC})$;
 i.\,e., $\varphi$ is true inside the Boolean valued universe
 ${\mathbb V}^{({\mathbb B})}$. \Endproc

 \smallskip

 We enter into the next agreement: If $\varphi (x)$ is a~formula of $\mathrm{ZFC}$ then,
 on assuming $x$ to be an element of
 ${\mathbb V}^{({\mathbb B})}$, the phrase ``$x$ satisfies $\varphi$ inside
 ${\mathbb V}^{({\mathbb B})}$'' or, briefly, ``$\varphi (x)$ is true inside
 ${\mathbb V}^{({\mathbb B})}$'' means that $[\![\varphi(x)]\!]=\mathbb 1$.
 This is sometimes written as ${\mathbb V}^{({\mathbb B})}\models \varphi (x)$.

 \subsubsec{1.1.4.} There is a~natural equivalence relation $x\sim y\iff[\![x=y]\!]=\mathbb 1$ in  the class
 ${\mathbb V}^{({\mathbb{B}})}$. Choosing a~representative of the least rank in each
 equivalence class or, more exactly, using the so-called ``Frege--Russell--Scott trick,''
 we obtain a~{\it separated Boolean valued universe\/}
 $\overline {\mathbb V}{}^{({\mathbb{B}})}$ for which
 $x=y\iff[\![x=y]\!]=\mathbb1$. It is easily to see that the Boolean truth value of a~formula
 remains unaltered if we replace in it each element of $\mathbb{V}^{({\mathbb B})}$ by one of its equivalents;
 cp.~\cite[\S~4.5]{IBA}. In this connection from now on we take ${\mathbb V}^{({\mathbb B})}\!:=
 \overline{\mathbb V}{}^{({\mathbb B})}$ without further specification.

 \subsubsec{1.1.5.} Given $x\in {\mathbb V}^{({\mathbb B})}$ and $b\in {\mathbb B}$, define the function
 $ bx:z\mapsto b\land x (z)$ $\bigl(z\in\dom(x)\bigr)$. Here we presume that $b\varnothing\!:=\varnothing$ for all
 $b\in {\mathbb B}$. Observe that in $\mathbb{V}^{({\mathbb B})}$ the  element $bx$ is defined correctly for
 $x\in\mathbb{V}^{({\mathbb B})}$ and $b\in {\mathbb B}$; cp.~\cite[\S~4.3]{IBA}.

 \subsubsec{1.1.6.} \proclaim{Mixing Principle.} Let
 $(b_{\xi})_{\xi\in\Xi}$ be a~{\it partition of unity\/} in~%
 ${\mathbb B}$, i.\,e., $\sup_{\xi\in\Xi} b_{\xi}=\mathbb 1$ and $\xi\neq\eta\Longrightarrow b_{\xi}\wedge
 b_{\eta}=\mathbb 0$. To each family $(x_{\xi})_{\xi\in\Xi}$
 in~$\mathbb{V}^{({\mathbb B})}$ there exists a~unique element $x$ in $\mathbb{V}^{({\mathbb B})}$
 such that $[\![x=x_{\xi}]\!]\ge b_{\xi}$ for all $\xi\in\Xi$.
 \Endproc

 This $x$ is  the {\it mixing\/} of
 $(x_{\xi})_{\xi\in\Xi}$ by~$(b_{\xi})_{\xi\in\Xi}$ denoted by $\mix_{\xi\in\Xi} b_{\xi} x_{\xi}$.

 \subsubsec{1.1.7.} \proclaim{Maximum Principle.} Let $\varphi(x)$ be
 a~formula of~$\mathrm{ZFC}$. Then {\rm(}in~$\mathrm{ZFC})$ there is a ${\mathbb B}$~valued set
 $x_0$ satisfying $[\![(\exists\,x)\varphi (x)]\!]=[\![\varphi
 (x_0)]\!]. $
 \Endproc

 In particular, if it is true within~$\mathbb{V}^{(\mathbb{B})}$
 that ``there is an $x$ for which $\varphi(x)$'', then there is an~element
 $x_0$ in~$\mathbb{V}^{(\mathbb{B})}$ (in the sense of $\mathbb{V}$)
 with $[\![\varphi(x_0)]\!]=\mathbb1$.
 In symbols, $\mathbb{V}^{(\mathbb{B})}\models
 (\exists\, x)\,\varphi(x)\,\Longrightarrow\, (\exists\, x_0)\,\mathbb{V}^{(\mathbb{B})}
 \models \varphi(x_0)$.

 \subsection{1.2.~Escher Rules}

 Now, we present a remarkable interplay between $\mathbb{V}$ and $\mathbb{V}^{(\mathbb{B})}$ which is based on the
 ope\-ra\-ti\-ons of canonical embedding, descent, and ascent.

 \subsubsec{1.2.1.}~We start with the canonical embedding of the von Neumann universe
 into the Boolean valued universe. Given ${x\in {\mathbb V}}$, we denote
 by~$x^{\scriptscriptstyle\wedge}$ the {\it standard name\/} of~$x$
 in~${\mathbb V}^{({\mathbb B})}$; i.e., the element defined by the
 following recursion schema:
 $$
 \varnothing^{\scriptscriptstyle\wedge}\!:=\varnothing,\quad \dom
 (x^{\scriptscriptstyle\wedge})\!:= \{y^{\scriptscriptstyle\wedge}:\,
 y\in x \}, \quad \operatorname{im}
 (x^{\scriptscriptstyle\wedge})\!:=\{\mathbb 1\}.
 $$
 Henceforth, working in the separated universe
 $\overline{\mathbb{V}}^{(\mathbb{B})}$, we agree to preserve the
 symbol $x^{\scriptscriptstyle\wedge}$ for the distinguished element
 of the class corresponding to~$x$. The map $x\mapsto x^{\scriptscriptstyle\wedge}$
 is called \textit{canonical embedding.}

 A~formula is  {\it bounded\/} or {\it restricted\/} provided that
 each bound variable in it is restricted by a~bounded quantifier;
 i.e., a~quantifier ranging over a~particular set. The latter means
 that each bound variable~$x$ is restricted by a~quantifier of the
 form~$(\forall\,x\in y)$ or $(\exists\,x\in y)$.

 \subsubsec{1.2.2.}
 \proclaim{Restricted Transfer Principle.}Let $\varphi(x_1,\dots,x_n)$
 be a bound\-ed formula of~$\mathrm{ZFC}$. Then $($in $\mathrm{ZFC})$ for every
 collection~$x_1,\dots,x_n\in\mathbb{V}$ we have
 $$
 \varphi(x_1,\dots,x_n)\,\Longleftrightarrow\,\mathbb{V}^{(\mathbb{B})}\models
 \varphi(x_1^{\scriptscriptstyle\wedge},\dots,x_n^{\scriptscriptstyle\wedge}).
 $$\Endproc

 \subsubsec{1.2.3.}~Given an arbitrary element $x$ of the Boolean valued universe
 $\mathbb{V}^{(\mathbb{B})}$, define the class $x{\downarrow}$ by
 $$
 x{\downarrow}\!:=\{y\in\mathbb{V}^{(\mathbb{B})}:\, [\![y\in x]\!]=\mathbb 1 \}.
 $$
 This class is called the {\it descent of\/} of~$x$. Moreover $x{\downarrow}$ is a set, i.\,e.,
 ${x{\downarrow}\in\mathbb{V}}$ for every~$x\in \mathbb{V}^{(\mathbb{B})}$. If $[\![x\ne\varnothing]\!]=\mathbb1$,
 then $x{\downarrow}$ is a non-empty set.

 \subsubsec{1.2.4.}~Suppose that $f$ is a map from~$X$ in~$Y$ within
 $\mathbb{V}^{(\mathbb{B})}$. More precisely, $f$, $X$ and $Y$ are in
 $\mathbb{V}^{(\mathbb{B})}$ and $[\![f:X\to Y]\!]=\mathbb1$. There exist a unique
 map $f{\downarrow}$ from $X{\downarrow}$ in~$Y{\downarrow}$
 (in the sense of the von Neumann universe $\mathbb{V}$) such that
 $$
 [\![f{\downarrow}(x)=f(x)]\!]=\mathbb 1\quad(x\in X{\downarrow}).
 $$
 Moreover, for a nonempty subset $A$ of $X$ within $\mathbb{V}^{(\mathbb{B})}$ (i.~e. $[\![\varnothing\ne A\subset X]\!]=\mathbb 1$) we have
 $f{\downarrow}(A{\downarrow})=f(A){\downarrow}$. The map
 $f{\downarrow}$ from $X{\downarrow}$
 to $Y{\downarrow}$ is called the {\it descent\/} of $f$ from~$\mathbb{V}^{(\mathbb{B})}$.
 The descent $f{\downarrow}$ of an~internal map $f$ is {\it extensional:}
 $$
 [\![x=x']\!]\leq[\![f{\downarrow}(x)=f{\downarrow}(x')]\!]
 \quad(x,x'\in X{\downarrow}).
 $$

 For the descents of the composite, inverse, and identity maps we have:
 $$
 (g\circ f){\downarrow}=g{\downarrow}\circ f{\downarrow},\quad
 (f^{-1}){\downarrow}=(f{\downarrow})^{-1},\quad
 (I_X){\downarrow}=I_{X{\downarrow}}.
 $$
 By virtue of~these rules we can consider the descent operation as a~functor from the category
 of~${\mathbb B}$-valued sets and mappings to the category of the standard sets and mappings
 (i.e., those in the sense of~${\mathbb V}$).

 \subsubsec{1.2.5.}~Given $x_1,\dots,x_n\in {\mathbb V}^{(\mathbb{B})}$, denote by
 $(x_1,\dots,x_n)^\mathbb{B}$ the corresponding
 ordered $n$-tuple inside~${\mathbb V}^{(\mathbb{B})}$. Assume
 that $P$ is an~$n$-ary relation on~$X$ inside~$\mathbb{V}^{({\mathbb B})}$; i.\,e.,
 $[\![P\subset X^{n^{\scriptscriptstyle\wedge}}]\!]=\mathbb1$
 and $[\![P\subset X^n]\!]=\nobreak{\mathbb 1}$.
 Then there exists an~$n$-ary relation~$P'$ on~$X{\downarrow}$ such
 that $(x_1,\dots,x_n)\in P'\iff [\![(x_1,\dots,x_n)^\mathbb{B}\in\nobreak P]\!]=\nobreak{\mathbb 1}$.
 Slightly abusing notation, we denote~$P'$ by the occupied symbol~$P{\downarrow}$ and
 call~$P{\downarrow}$ the {\it descent\/} of~$P$.

 \subsubsec{1.2.6.}~Let $x\in {\mathbb V}$ and $x\subset {\mathbb
 V}^{(\mathbb{B})}$; i.\,e., let $x$ be some set composed of $\mathbb{B}$-valued
 sets or, symbolically, $x\in \mathcal{P}({\mathbb V}^{(\mathbb{B})})$.
 Put $\varnothing{\uparrow}\!:=\varnothing$ and $\dom(x{\uparrow})\!:=x$,
 $\im (x{\uparrow})\!:=\{{\mathbb 1}\} $ if $x\neq \varnothing$.
 The element~$x{\uparrow}$ \big(of the non-separated
 universe~$\overline{\mathbb{V}}^{({\mathbb B})}$, i.\,e., the distinguished
 representative of the class~$\{y\in\overline{\mathbb{V}}^{({\mathbb B})}:
 [\![y=\nobreak x{\uparrow}]\!]=\nobreak\mathbb{1}\}$\big)
 is  the {\it ascent\/} of~$x$.
 For the corresponding element in the separated universe
 $\mathbb{V}^{(\mathbb{B})}$ the same name and notation are preserved.

 \subsubsec{1.2.7.}~Let $X,Y,f\in\mathcal{P}(\mathbb{V}^{(\mathbb{B})})$ and
 $f$ be a mapping from~$X$ to~$Y$. There exists a mapping
 $f{\uparrow}$ from $X{\uparrow}$ to $Y{\uparrow}$ within $\mathbb{V}^{(\mathbb{B})}$
 satisfying
 $$
 [\![f{\uparrow}(x)=f(x)]\!]=\mathbb{1}\quad(x\in X),
 $$
 if and only if $f$ is {\it extensional}, i.\,e., the relation
 holds:
 $$
 [\![x=x']\!]\leq [\![f(x)=f(x')]\!]\quad(x,x'\in X).
 $$
 The map $f{\uparrow}$ with the above property is unique and satisfy the relation
 $f{\uparrow}(A{\uparrow})=f(A){\uparrow}$ $(A\subset X)$.
 The composite of extensional maps is extensional. Moreover, the ascent of a~composite
 is equal to the composite of the ascents inside~$\mathbb{V}^{(\mathbb{B})}$:
 $$
 \mathbb{V}^{(\mathbb{B})}\models(g\circ f){\uparrow}=g{\uparrow}\circ f{\uparrow}.
 $$
 Observe also that if $f$ and $f^{-1}$ are extensional then
 $(f{\uparrow})^{-1}=(f^{-1}){\uparrow}$.

 \subsubsec{1.2.8.}~Suppose that $X\in \mathbb{V}$, $X\ne\varnothing$; i.\,e.,
 $X$ is a~nonempty set. Let the letter~$\iota$ denote the standard name embedding
 $x\mapsto x^{\scriptscriptstyle\wedge}$ $(x\in X)$.
 Then $\iota(X){\uparrow}=X^{\scriptscriptstyle\wedge}$ and
 $X=\iota^{-1}(X^{\scriptscriptstyle\wedge}{\downarrow})$. Take $Y\in\mathbb{V}^{(\mathbb{B})}$
 with $[\![Y\ne\varnothing]\!]=\mathbb{1}$. Using the
 above relations, we may extend the ascent operation to the case of maps from $X$ to $Y{\downarrow}$
 and descent operation to
 the case of internal maps from~$X^{\scriptscriptstyle\wedge}$
 to~$Y]\!]=\mathbb{1}$.

 The maps $f\upwardarrow$ and $g\downwardarrow$ are called {\it modified
 descent\/} of $f$ and {\it modified ascent\/} of~$g$, respectively.
 (Again, when there is no ambiguity, we simply speak of ascents and use simple
 arrows.) It is easy to see that $g\upwardarrow$ is the unique map from $X$
 to $Y{\downarrow}$ satisfying
 $$
 g\downwardarrow(x)=g(x^{\scriptscriptstyle\wedge}){\downarrow}
 \quad (x\in X)
 $$
 and $f{\upwardarrow}$ is the~unique map
 from $X^{\scriptscriptstyle\wedge}$ to $Y$ within $\mathbb{V}^{(\mathbb{B})}$
 satisfying
 $$
 [\![f\upwardarrow(x^{\scriptscriptstyle\wedge})=f(x)]\!]=\mathbb{1} \quad(x\in X).
 $$

 \subsubsec{1.2.9.}~Given  $X\subset \mathbb{V}^{(\mathbb{B})}$, we denote by~$\operatorname{mix}(X)$ the set of all mixtures of the form~$\operatorname{mix} (b_{\xi}x_{\xi})$, where $(x_{\xi})\subset X$
 and $(b_{\xi})$ is an~arbitrary partition of unity. The following
 assertions are referred to as the {\it rules for canceling
 arrows\/} or the {\it Escher rules}.

 Let $X$ and $X'$ be subsets of~$\mathbb{V}^{(\mathbb{B})}$ and let
 $f:X\to X'$ be an~extensional mapping. Suppose that $Y,Y',g\in
 \mathbb{V}^{(\mathbb{B})}$ are such that $[\![\,Y\ne \varnothing]\!]=[\![\,g:Y\to
 Y']\!]=\mathbb{1}$. Then
 $$
 \gathered
 X{\uparrow}{\downarrow}=\operatorname{mix}(X),\quad
 Y{\downarrow}{\uparrow}=Y;\\
 f{\uparrow}{\downarrow}=f,\quad g{\downarrow}{\uparrow}=g.
 \endgathered
 $$
 There are some other cancelation rules.

 \subsection{1.3. Boolean Valued Reals and Vector Lattices}

 The main results of the section tells us that the Boolean valued interpretation of the field of reals (complexes)
 is a real (complex) universally complete vector lattice.
 Everywhere below $\mathbb{B}$ is a complete Boolean algebra and
 ${\mathbb{V}}^{(\mathbb{B})}$ is the corresponding
 Boolean valued universe.

 \subsubsec{1.3.1.}~By virtue of the Transfer and Maximum Principles there exists
 an~element~$\mathrsfs R\in\mathbb{V}^{(\mathbb{B})}$ for which $[\![\,\mathrsfs R$ is a~field of
 reals$\,]\!]=\mathbb{1}$. Note also that $\varphi (x)$, formally presenting the expressions of
 the axioms of an~archimedean ordered field $x$, is bounded;
 therefore, by the Restricted Transfer Principle $[\![\,\varphi (\mathbb R^{\land})\,]\!]=\mathbb 1$, i.\,e.,
 $[\![\,\mathbb{R}^{\land}$ is an~archimedean ordered
 field$\,]\!]=\mathbb{1}$. Thus, we will assume that $\mathbb{R}^{\scriptscriptstyle\wedge}$ is
 a~dense subfield of $\mathrsfs R$, while the elements~$0\!:=0^{\scriptscriptstyle\wedge}$
 and~$1\!:=1^{\scriptscriptstyle\wedge}$ are the zero and unity of~$\mathrsfs R$ within the
 model~$\mathbb{V}^{(\mathbb{B})}$.

 \subsubsec{1.3.2.}~Let $\circledR$ be the underlying set of the field $\mathcal{R}$ on which the addition
 $\oplus$, multiplication $\otimes$, and ordering $\oleq$ are given.
 Then $\mathcal{R}$ is a $5$-tuple $(\circledR,\oplus, \otimes, \oleq,
 0^{\scriptscriptstyle\wedge},1^{\scriptscriptstyle\wedge})$ within $\mathbb{V}^{(\mathbb{B})}$; in symbols,
 $\mathbb{V}^{(\mathbb{B})}\models\mathcal{R}=(\circledR,\oplus, \otimes, \oleq, 0^{\scriptscriptstyle\wedge}, 1^{\scriptscriptstyle\wedge})$.

 The descent $\mathcal{R}{\downarrow}$ of the field
 $\mathcal{R}$ is the descent of the underlying set $\circledR{\downarrow}$
 together with the descended operations $+\!:=\oplus{\downarrow}$,
 $\cdot\!:=\odot{\downarrow}$, order relation $\preccurlyeq\!:=\oleq{\downarrow}$,
 and distinguished elements
 $0^{\scriptscriptstyle\wedge},\,1^{\scriptscriptstyle\wedge}$;
 in symbols, $\mathcal{R}{\downarrow}=(\circledR,\mathbf{R}{\downarrow},\oplus{\downarrow},
 \otimes{\downarrow}, \oleq {\downarrow},
 0^{\scriptscriptstyle\wedge},1^{\scriptscriptstyle\wedge})$.
 Also, we may introduce multiplication by the standard reals in~$\mathrsfs R{\downarrow}$ by the rule
 $$
 y=\lambda x\Longleftrightarrow
 [\![\,y=\lambda^{\scriptscriptstyle\wedge}\odot x\,]\!]=\mathbb{1}\quad
 (\lambda\in \mathbb R,\ x, y\in|\mathcal{R}|{\downarrow}).
 $$

 \subsubsec{1.3.3.} \proclaim{Gordon Theorem.}Let $\mathrsfs R$
 be the reals in~$\mathbb{V}^{(\mathbb{B})}$. Then $\mathrsfs
 R{\downarrow}$ $($with the descended operations and order$)$ is
 a~universally complete vector lattice with a weak order unit~$\mathbb{1}\!:=1^{\scriptscriptstyle\wedge}$.
 Moreover, there exists a Boolean isomorphism~$\chi$ of~$\mathbb{B}$ onto
 the base~$\mathbb{P}(\mathrsfs R{\downarrow})$ such that for all
 $x, y\in \mathrsfs R{\downarrow}$ and $b\in\mathbb{B}$ we have
  $$
 \gathered \chi (b) x=\chi (b) y\Longleftrightarrow b\le [\![\,x=y\,]\!],
 \\
 \chi (b) x\le \chi (b) y\Longleftrightarrow b\le [\![\,x\le
 y\,]\!].
 \endgathered
 \eqno{(G)}
 $$
 \Endproc

 \subsubsec{1.3.4.}~A~vector lattice is an~{\it $f$-algebra\/} if it is simultaneously a real algebra and satisfies, for all $a, x, y\in X_+$, the conditions: 1)~$x\geq0$ and $y\geq0$ imply $xy\geq0$ and 2)~$x\perp y=0$ implies that $(ax) \perp y$ and $(xa)\perp y$. The multiplication in every (Archimedean) $f$-algebra is commutative and associative.  An $f$-algebra is called {\it semi-prime\/} if $xy=0$ implies $x\perp y$ for all  $x$ and $y$.
 The universally complete vector lattice $\mathcal{R}{\downarrow}$ with the descended multiplications is a~semiprime~$f$-algebra with a ring unit~$\mathbb 1\!=1^{\scriptscriptstyle\wedge}$.

 \subsubsec{1.3.5.}~By the Maximum Principle 1.1.7, there is an~element $\mathcal C\in {\mathbb V}^{(\mathbb B)}$
 for which $[\![\,\mathcal C$ is the complexes$\,]\!]
 ={\mathbb{1}}$. Since the equality
 ${\mathbb C}={\mathbb R}\oplus i{\mathbb R}$ is expressed
 by a~bounded set-theoretic formula, from the Restricted
 Transfer Principle 1.2.2 we obtain $[\![\,{\mathbb C}^{\scriptscriptstyle\wedge} =
 {\mathbb R}^{\scriptscriptstyle\wedge}\oplus
 i^{\scriptscriptstyle\wedge}{\mathbb R}^{\scriptscriptstyle\wedge}\,]\!]=
 {\mathbb{1}}$. Moreover, ${\mathbb R}^{\scriptscriptstyle\wedge}$
 is assumed to be a~dense subfield of~$\mathcal R$;
 therefore, we can also assume that
 ${\mathbb C}^{\scriptscriptstyle\wedge}$ is a~dense subfield
 of~$\mathcal C$. If $1$ is the unity of~$\mathbb C$ then
 $1^{\scriptscriptstyle\wedge}$ is the unity of~$\mathcal
 C$ inside~$\mathbb V^{(\mathbb B)}$. We write $i$ instead
 of~$i^{\scriptscriptstyle\wedge}$ and $\mathbb{1}$ instead of
 $1^{\scriptscriptstyle\wedge}$.
 By the Gordon Theorem $\mathcal C{\downarrow}=\mathcal
 R{\downarrow}\oplus i\mathcal R{\downarrow}$; consequently
 $\mathcal C{\downarrow}$ is a~universally complete complex
 vector lattice, i.~e., the complexification of a vector lattice $\mathcal{R}{\downarrow}$.
 Moreover, $\mathcal C{\downarrow}$ is a~complex $f$-algebra
 defined as the complexification of a real $f$-algebra
 with a~ring unit $\mathbb{1}\!:=1^{\scriptscriptstyle\wedge}$.

 \subsubsec{1.3.6.}~Let $A$ be an $f$-algebra. A vector lattice
 $X$ is said to be an $f$-\textit{module over\/} $A$ if the following hold:
 \vspace{2pt}

 \subsubsec{(1)}~$X$ is a module over $A$ (with respect to a multiplication $A\times X\ni(a,x)\mapsto ax\in X$);
 \vspace{2pt}

 \subsubsec{(2)}~$ax\geq 0$ for all $a\in A_+$ and $x\in X_+$;
 \vspace{2pt}

 \subsubsec{(3)}~$x\perp y$ implies $ax\perp y$ for all $a\in A_+$
 and $x,y\in X$.
 \vspace{2pt}

 A vector lattice $X$ has a natural $f$-module structure  over
 $\Orth(X)$, i.~e., $\pi x\!:=\pi(x)$ for all $x\in X$ and $\pi\in\Orth(X))$. Clearly,
 $X$ is an $f$-module over an arbitrary $f$-sub\-mo\-dule $A\subset\Orth(X)$ and, in particular, over $\mathcal{Z}(X)$.

 \subsubsec{1.3.7.}~\proclaim{Theorem}~Let $X$ be an $f$-module over $\mathcal{Z}(Y)$ with
 $Y$ being a Dedekind complete vector lattice and $\mathbb{B}=\mathbb{P}(Y)$. Then there exists
 $\mathcal{X}\in\mathbb{V}^{(\mathbb{B})}$ such that $[\![\mathcal{X}$~is a vector lattice over
 $\mathcal{R}]\!]=\mathbb{1}$, $\mathcal{X}{\downarrow}$ is an
 $f$-module over $A^{\mathrm{u}}$, and there is an $f$-module isomorphism $h$ from $X$ to
 $\mathcal{X}{\downarrow}$ satisfying $\mathcal{X}{\downarrow}=\mix(h(X))$.
 \Endproc

 \subsubsec{1.3.8.}~Theorem 1.3.3 was established in \cite{Gor1}.
 The concept of an $f$-module was introduced in \cite{LS}.

 \subsection{1.4. Boolean Valued Functionals}

 We will demonstrate in this section how Boolean valued analysis works by trans\-fer\-ring
 some results from order bounded functionals to operators. Below $X$ and $Y$ stand for vector
 lattices, where~$Y$ is an order dense sublattice in $\mathcal{R}{\downarrow}$.

 \subsubsec{1.4.1.}~Let $\mathbb{B}$ be a complete Boolean algebra and
 $\mathcal{R}$ be the field of reals in $\mathbb{V}^{(\mathbb{B})}$.
 The fact that $X$ is a~vector lattice over the ordered
 field $\mathbb{R}$ may be rewritten as a~restricted formula, say, $\varphi (X, \mathbb{R})$.
 Hence, recalling the Restricted Transfer Principle~1.2.2, we come to the identity
 $[\![\,\varphi (X^{\scriptscriptstyle\wedge},
 \mathbb{R}^{\scriptscriptstyle\wedge})\,]\!]=\mathbb{1}$ which amounts to saying that
 $X^{\scriptscriptstyle\wedge}$ is a~vector lattice over the ordered
 field~$\mathbb{R}^{\scriptscriptstyle\wedge}$ inside $\mathbb{V}^{(\mathbb{B})}$.

 Let $X^{{\scriptscriptstyle\wedge}\sim}\!:= L^{\sim}
 (X^{\scriptscriptstyle\wedge},\mathcal{R})$~be the space
 of order bounded $\mathbb{R}^{\scriptscriptstyle\wedge}$-linear
 functionals from $X^{\scriptscriptstyle\wedge}$ to
 $\mathcal{R}$. More precisely, $\mathcal{R}$ is considered
 as a vector space over the field
 $\mathbb{R}^{\scriptscriptstyle\wedge}$ and by Maximum Principle there exists
 ${X^{\scriptscriptstyle\wedge\sim}\in\mathrm{V}^{(\mathbb{B})}}$
 such that $[\![X^{\scriptscriptstyle\wedge\sim}$~is
 a vector space over $\mathcal{R}$ of $\mathbb{R}^{\scriptscriptstyle\wedge}$-linear
 order bounded functionals from $X^{\scriptscriptstyle\wedge}$
 to~$\mathcal{R}$ ordered by the cone of positive
 functionals~$]\!]=\mathbb{1}$. A functional $\tau\in X^{\scriptscriptstyle\wedge\sim}$ is positive if
 $[\![(\forall\,x\in X^{\scriptscriptstyle\wedge}) \tau(x)\geq0]\!]=\mathbb{1}$.

 \subsubsec{1.4.2.}~It can be easily seen that the Riesz--Kantorovich Theorem remains
 true if $X$ is a vector lattice over a dense subfield $\mathbb{P}\subset\mathbb{R}$ and
 $Y$ be a Dedekind complete vector lattices (over $\mathbb{R}$) and $L^{\scriptscriptstyle\sim}(X,Y)$
 is replaced by $L_{\mathbb{P}}^{\scriptscriptstyle\sim}(X,Y)$, the vector spaces over
 $\mathbb{R}$ of all $\mathbb{P}$-linear order bounded operators from $X$ to $Y$, ordered
 by the cone of positive operators: $L_{\mathbb{P}}^{\scriptscriptstyle\sim}(X,Y)$ is a
 Dedekind complete vector lattice. Thus, the descent $X^{\scriptscriptstyle\wedge\sim}{\downarrow}$
 of $X^{\scriptscriptstyle\wedge\sim}$ is a~Dedekind complete
 vector lattice. Boolean valued interpretation of this fact yields that $X^{\scriptscriptstyle\wedge}\!:=
 L_{\mathbb{R}^{\scriptscriptstyle\wedge}}^\sim(X^{\scriptscriptstyle\wedge},\mathbb{R})$
 is a Dedekind complete vector lattice within $\mathbb{V}^{(\mathbb{B})}$ with $\mathbb{B}\!:=\mathbb{P}(Y)$.
 In particular, the descent $X^{{\scriptscriptstyle\wedge}\sim}{\downarrow}$ of the space
 $X^{{\scriptscriptstyle\wedge}\sim}$ is Dedekind complete vector lattice. Let $L_{dp}^\sim(X,Y)$
 and $\Hom(X,T)$ stand respectively for the space of disjointness preserving order bounded
 operators and the set of all lattice homomorphism   from $X$ to $Y$. The following result
 is based on the construction from~1.2.8.

 \subsubsec{1.4.3.}~\theorem{}Let $X$ and $Y$ be vector lattices
 with $Y$ universally complete and represented as
 $Y=\mathcal{R}{\downarrow}$. Given $T\in L^{\sim}(X,Y)$,
 the modified ascent $T{\upwardarrow}$ is an~order bounded
 $\mathbb{R}^{\scriptscriptstyle\wedge}$-linear functional on
 $X^{\scriptscriptstyle\wedge}$ within
 $\mathrm{V}^{(\mathbb{B})}$; i.~e.,
 $[\![\,T{\upwardarrow}\in
 X^{\scriptscriptstyle\wedge\sim}\,]\!]
 =\mathbb{1}$. The mapping $T\mapsto T{\upwardarrow}$ is a
 lattice isomorphism between the Dedekind complete vector
 lattices $L^{\sim}(X,Y)$ and
 $X^{\scriptscriptstyle\wedge\sim}{\downarrow}$.
 \Endproc

 \subsubsec{1.4.4.}~\proclaim{Corollary.}Given operators $R,S\in L^{\scriptscriptstyle\sim}(X, Y)$,
 put $\sigma\!:=S{\upwardarrow}$ and $\tau\!:=T{\upwardarrow}$. The following are hold true:
 \vspace{2pt}

 \subsubsec{(1)}~$S\leq T\,\Longleftrightarrow\,[\![\,\sigma\leq\tau\,]\!]
 =\mathbb{1};$
 \vspace{2pt}

 \subsubsec{(2)}~$S=|T|\,\Longleftrightarrow\,[\![\,\sigma=|\tau|\,]\!]
 =\mathbb{1};$
 \vspace{2pt}

 \subsubsec{(3)}~$S\perp T\,\Longleftrightarrow\,[\![\,\sigma\perp\tau\,]\!]
 =\mathbb{1};$
 \vspace{2pt}

 \subsubsec{(4)}~$[\![\,T\in\Hom(X,Y)\,]\!]\,\Longleftrightarrow\,
 [\![\,\tau\in\Hom(X^{\scriptscriptstyle\wedge},\mathcal{R})\,]\!]=\mathbb{1}$;
 \vspace{2pt}

 \subsubsec{(5)}~$T \in L_{dp}^{{\scriptscriptstyle}\sim} (X,Y)\,\Longleftrightarrow\,
 [\![\,\tau \in(X^{\scriptscriptstyle\wedge\sim})_{dp}\,]\!]=\mathbb{1}$.
 \Endproc

 \subsubsec{1.4.5.}~Consider a vector lattice $X$ and let
 $D$ be an order ideal in $X$. A~linear operator $T$ from $D$ into $X$ is
 \textit{band preserving} provided that one (and hence all) of the following  holds:
 $x\perp y$ implies $Tx\perp y$ for all $x\in D$ and $y\in X$, or
 equivalently, $Tx\in\{x\}^{\perp\perp}$ for all $x\in D$ (the disjoint complements
 are taken in $X$). If $X$ is a vector lattice with the principal projection
 property and $D\subset X$ is an order dense ideal, then
 a linear operator $T:D\to X$ is band preserving if and
 only if $T$ commutes with  band projections:
 $\pi Tx=T\pi x$ for all $\pi\in\mathbb{P}(X)$ and $x\in D$.

 \subsubsec{1.4.6.}~Let $\End_N(X_{\mathbb{C}})$~be the set of all band preserving endomorphisms
 of~$X_{\mathbb{C}}$ with $X\!:=\mathcal{R}{\downarrow}$. Clearly,
 $\End_N(XG_{\mathbb{C}})$ is a~complex vector space. Moreover,
 $\End_N(X_{\mathbb{C}})$ becomes a faithful unitary module over the ring $X_{\mathbb{C}}$ on letting
 $gT$ be equal to $gT:x\mapsto g\cdot Tx$ for all $x\in X_{\mathbb{C}}$. This is
 immediate since the multiplication by an element of $X_{\mathbb{C}}$ is band
 preserving and the composite of band preserving operators is band
 preserving too.

 \subsubsec{1.4.7.}~By $\End_{\mathbb{C}^{\scriptscriptstyle\wedge}}(\mathcal{C})$
 we denote the element of $\mathbb{V}^{(\mathbb{B})}$ that represents
 the space of all $\mathbb{C}^{\scriptscriptstyle\wedge}$-linear
 operators from  $\mathcal{C}$ into $\mathcal{C}$. Then
 $\End_{\mathbb{C}^{\scriptscriptstyle\wedge}}(\mathcal{C})$ is a~vector
 space over $\mathcal{C}$ inside $\mathbb{V}^{(\mathbb{B})}$, and
 $\End_{\mathbb{C}^{\scriptscriptstyle\wedge}}(\mathcal{C}){\downarrow}$ is
 a~faithful unitary module over a complex $f$-slgebra~$X_\mathbb{C}$.

 \subsubsec{1.4.8.}~\proclaim{}A linear operator~$T$ on a universally
 complete vector lattice $X$ or $X_\mathbb{C}$ is band pre\-ser\-ving if
 and only if $T$ is extensional.
 \Endproc

 \beginproof~Take a linear operator $T:X\to X$. By the Gordon Theorem the ex\-ten\-sionality condition
 $[\![x=y]\!]\leq[\![Tx=Ty]\!]$ $(x,y\in X=\mathcal{R}{\downarrow})$ amounts to saying that
 the identity $\pi x=\pi y$ implies $\pi Tx=\pi Ty$ for all
 $x,y\in X$ and $\pi \in\mathbb{P}(X)$. By linearity of $T$ the latter is equivalent to
 $\pi x=0\,\Longrightarrow\,\pi Tx=0$ $(x\in X$, $\pi\in \mathbb{P}(X))$.
 Substituting $y\!:=\pi^\perp y$ yields $\pi T\pi^\perp=0$ or, which is the same,
 $\pi T=\pi T\pi$. According to 1.4.5 $T$ is band preserving. The complex case is
 treated by complexification.~\endproof

 \subsubsec{1.4.9.}~\proclaim{Theoem.}The modules $\End_N(X_\mathbb{C})$ and $\End_{\mathbb{C}^{\scriptscriptstyle\wedge}}(\mathcal{C}){\downarrow}$
 are isomorphic. The isomorphy may be established by sending a band preserving operator to its ascent.
 The same remains true when $\mathcal{C}$ and $\mathbb{C}$ are replaced by $\mathcal{R}$ and $\mathbb{R}$, respectively.
 \Endproc

 \beginproof~By virtue of 1.4.8, we can apply the constructions 1.2.4 and 1.2.7, as well as the cancelation
 rules~1.2.9.~\endproof

 \section{Chapter 2. Band Preserving Operators}

 \subsection{2.1. Wickstead's Problem and Cauchy's Functional\\ Equation}

 In this section we demonstrate that the band preserving operators in universally complete
 vector lattices are solutions in disguise of the Cauchy functional equation and the Wickstead
 problem amounts to that of regularity of all solutions to the equation.

 \subsubsec{2.1.1.}~\proclaim{The Wickstead Problem:}When are we so happy in a vector lattice
 that all band preserving linear operators turn out to be order bounded?
 \Endproc

 This question was raised by~Wickstead in~\cite{Wic1}. Further progress is presented in~\cite{AVK, AVK1,  Gut6, Kus11, Kus2, MW}. Combined approach based on logical, algebraic, and analytic tools was presented in \cite{Kus11, Kus2, Kus4}.
 A survey of the main ideas and results on the problem and its modifications see in~\cite{GKK}.

 The answer depends on the~vector lattice in which the~operator in question acts. Therefore,
 the problem can be reformulated as follows: \textsl{Characterize the vector lattices
 in which every band preserving linear operators is order bounded.}

 Let $X$ be a universally complete vector lattice and $T$ a band preserving linear operator in~$X$.
 By the Gordon Theorem we may assume that $X=\mathcal{R}{\downarrow}$, where $\mathcal{R}$
 is the field of reals within $\mathbb{V}^{(\mathbb{B})}$ and
 $\mathbb{B}=\mathbb{P}(X)$. Moreover, according to Theorem~1.4.9 we can assume further
 that $T=\tau{\downarrow}$, where $\tau\in\mathbb{V}^{(\mathbb{B})}$ is an internal
 $\mathbb{R}^{\scriptscriptstyle\wedge}$-linear function from  $\mathcal{R}$
 to $\mathcal{R}$. It can be easily seen that $T$ is order bounded if and only if $[\![\tau$
 is order bounded (i.\,e., $\tau$ is bounded on intervals
 $[a,b]\subset\mathcal{R})]\!]=\mathbb{1}$.

 \subsubsec{2.1.2.}~By $\mathbb{F}$ we denote either $\mathbb{R}$ or $\mathbb{C}$.
 The {\it Cauchy functional equation\/} with
  $f:\,\mathbb{F}\to\mathbb{F}$ unknown has the form
 $$
 f(x+y)=f(x)+f(y)\quad (x,y\in\mathbb{F}).
 $$
 It is easy that a solution to the equation is automatically $\mathbb{Q}$-homogeneous, i.\,e. it satisfies
 another functional equation:
 $$
 f(qx)=qf(x)\quad (q\in\mathbb{Q},\ x\in\mathbb{F}).
 $$
 In the sequel we will be interested in a more general situation. Namely, we will
 consider the simultaneous functional equations
 $$
 \begin{cases}
 f(x+y)=f(x)+f(y)\quad (x,y\in\mathbb{F}),\\
 f(px)=pf(x)\quad (p\in\mathbb{P},\ x\in\mathbb{F}),
 \end{cases}
 \eqno(L)
 $$
 where $\mathbb{P}$ is a subfield of~$\mathbb{F}$ which includes~$\mathbb{Q}$.
 Denote by $\mathbb{F}_{\mathbb{P}}$ the field
 $\mathbb{F}$ which is considered as a vector space over~$\mathbb{P}$.
 Clearly, solutions to the simultaneous equations $(L)$ are precisely
  $\mathbb{P}$-linear functions from~$\mathbb{F}_{\mathbb{P}}$ to~$\mathbb{F}_{\mathbb{P}}$.

 \subsubsec{2.1.3.}~Let $\mathcal{E}$ be a~Hamel basis for a~vector space~$\mathbb{F}_{\mathbb{P}}$,
 and let $\mathcal{F}(\mathcal{E},\mathbb{F})$ be the~space of all functions from~$\mathcal{E}$ to~$\mathbb{F}$.
 The solution set of~$(L)$ is a~vector space over~$\mathbb{F}$ isomorphic with~$\mathcal{F}(\mathcal{E},\mathbb{F})$.
 Such an isomorphism can be implemented by sending a solution~$f$ to the restriction $f|_{\mathcal{E}}$ of~$f$~to~$\mathcal{E}$.
 The inverse
 isomorphism $\varphi\mapsto f_\varphi$ $\big(\varphi\in\mathcal{F}(\mathcal{E},\mathbb{F})\big)$ is defined by
 $$
 f_{\varphi}(x):=\sum\limits_{e\in\mathcal{E}}\varphi(e)\psi(e)
 \quad (x\in\mathbb{F}_{\mathbb{P}}),
 $$
 where $x=\sum_{e\in\mathcal{E}}\psi(e)e$ is the expansion of $x$ with respect to Hamel basis $\mathcal{E}$.

 \subsubsec{2.1.4.} \proclaim{Theorem}Each solution of~$(L)$ is either $\mathbb{F}$-linear or everywhere
 dense in~$\mathbb{F}^2:=\mathbb{F}\times\mathbb{F}$. In particular if $\mathcal{E}$ is a Hamel basis for
 a~vector space~$\mathbb{F}_{\mathbb{P}}$ and $f$ be the unique extension of a function $\varphi:\mathcal{E}\to\mathbb{R}$. Then
 $f$~is continuous if and only if $\varphi(e)/e={\rm const}$ $(e\in\mathcal{E})$.
 \Endproc

 \subsubsec{2.1.5.} Assume now that $\mathbb{F=\mathbb{C}}$ and $\mathbb{P}:=\mathbb{P}_0+i\mathbb{P}_0$
 with $\mathbb{P}_0$ being a subfield in~$\mathbb{R}$. Then the space of solutions of the system $(L)$ is
 a complexification of the space of solution of the same system with $\mathbb{P}:=\mathbb{P}_0$.
 In more detail, if $g:\,\mathbb{R}\to\mathbb{R}$ is a $\mathbb{P}_0$-linear function then we have the unique $\mathbb{P}$-linear function $\widetilde{g}:\,\mathbb{C}\to\mathbb{C}$ defined as
 $$
 \widetilde{g}(z)=g(x)+ig(y)\quad(z=x+iy\in\mathbb{C}).
 $$
 Conversely, if $f:\,\mathbb{C}\to\mathbb{C}$ is a $\mathbb{P}$-linear function
 then we have the a unique pair of $\mathbb{P}_0$-linear functions $g_1,g_2:\,\mathbb{R}\rightarrow \mathbb{R}$
 such that $f(z)=\widetilde{g}_1(z)+i\widetilde{g}_2(z)$ $(z\in\mathbb{C})$. Thus, every solution $f$ of $(L)$
 can be represented in the form $f=f_1+if_2$, where $f_1,f_2:\,\mathbb{C}\to\mathbb{C}$ are $\mathbb{P}_0$-linear
 and $f_i(\mathbb{R})\subset\mathbb{R}$ $(i=1,2)$ have the same property. Say that $f$ is {\it monotone\/} or {\it bounded\/} if $f_1$ and $f_2$.

 \subsubsec{2.1.6.} \proclaim{}Let $\mathbb{P}$ be a subfield of~$\mathbb{F}$, while $\mathbb{P}:=\mathbb{P}_0+i\mathbb{P}_0$
 for some dense subfield $\mathbb{P}_0\subset\mathbb{R}$, in case $\mathbb{F}=\mathbb{C}$. The following are equivalent:

 \subsubsec{(1)} $\mathbb{F}=\mathbb{P}$;

 \smallskip

 \subsubsec{(2)} every solution to~$(L)$ is order bounded.\Endproc

 \beginproof~The implication  $(1)\Longrightarrow(2)$ is trivial. Prove the converse by way of contradiction.
 The assumption that $\mathbb{F}\ne\mathbb{P}$ implies that each Hamel basis $\mathcal{E}$ for the vector space
 $\mathbb{F}_{\mathbb{P}}$ contains at least two nonzero distinct elements $e_1,e_2\in\mathcal{E}$.
 Define the function $\psi:\,\mathcal{E}\to\mathbb{F}$ so that $\psi(e_1)/e_1\ne\psi(e_2)/e_2$. Then
 the $\mathbb{P}$-linear function $f=f_{\psi}:\,\mathbb{F}\to\mathbb{F}$, coinciding with~$\psi$
 on~$\mathcal{E}$, would exist by~2.1.2 and be discontinuous by~2.1.4.~\endproof

 \subsubsec{2.1.7.} Add to the system $(L)$ the equation $f(xy)=f(x)f(y)$ (or $f(xy)=f(x)y+xf(y)$)
 $(x,y\in\mathbb{F})$. A solution of the resulting system is called  $\mathbb{P}$-\textit{endomorphism\/}
 ($\mathbb{P}$-\textit{derivation\/}). The existence of the nontrivial $\mathbb{P}$-endomorphism and
 $\mathbb{P}$-derivation can be obtained similarly, but using a transcendental basis instead of a Hamel
 basis (cp.~\cite{Kuc}). Interpreting such existence results in a Boolean valued model yields the
 existence of band preserving endomorphism and  derivations of a universally
complete $f$-algebra,
 see \cite{Kus2, Kus4}, as well as \cite{IBA}.

 \section{2.2. Locally One-Dimensional Vector Lattices}

 Boolean valued representation of a vector lattice is a vector sublattice in $\mathcal{R}$ considered as a vector lattice over $\mathbb{R}^{\scriptscriptstyle\wedge}$. It stands to reason to find out what construction in a vector lattice corresponds to a Hamel basis for its Boolean valued representation.

 \subsubsec{2.2.1.}~Let $X$ be a vector lattice with a cofinal
 family of band projections. We will say that  $x,y\in X$ {\it differ\/} at $\pi\in\mathbb{P}(X)$
 provided that $\pi|x-y|$ is a weak order unit in~$\pi(X)$  or, equivalently, if  $\pi (X)\subset|x-y|^{\perp\perp}$.
 Clearly, $x$ and $y$ differ at $\pi$ whenever $\rho x=\rho y$ implies $\pi\rho=0$
 for all $\rho\in\mathbb{P}(X)$.
 A subset $\mathcal E$ of~$X$ is said to be {\it locally linearly independent\/}
 provided that, for an arbitrary nonzero band projection
 $\pi$  in~$X$ and each collection of the elements $e_1,\ldots,e_n\in\mathcal E$
 that are pairwise different at $\pi$, and each collection of reals
 $\lambda_1,\ldots,\lambda_n\in\mathbb R$, the condition
 $\pi(\lambda_1e_1+\cdots+\lambda_ne_n)=0$ implies that
 $\lambda_k=0$ for all $k\!:=1,\ldots,n$. In other words,
 $\mathcal E$  is locally linearly independent if for every band
 projection $\pi\in\mathbb{P}(X)$ any subset of $\pi(\mathcal{E})$
 consisting of nonzero members pairwise different at $\pi$ is linearly independent.

 An inclusion-maximal
 locally linearly independent subset of~$X$ is called a~{\it local Hamel basis\/} for~$X$.

 \subsubsec{2.2.2.}~\proclaim{} Each  vector lattice $X$ with a cofinal family of band
 projections has a local Hamel basis for~$X$.
 \Endproc

 \subsubsec{2.2.3.}~A~locally linearly independent set
 $\mathcal E$ in~$G$ is a~local Hamel basis if and only if for every $x\in G$ there exist
 a partition of unity $(\pi_\xi)_{\xi\in\Xi}$ in $\mathbb{P}(G)$ and a family of reals
 $(\lambda_{\xi,e})_{\xi\in\Xi,e\in\mathcal E}$ such that
 $$
 x=\osum_{\xi\in\Xi}\left(\sum_{e\in\mathcal E}
 \lambda _{\xi ,e}\pi_\xi e\right)
 $$
 and for every $\xi\in\Xi$ the set $\{e\in\mathcal E:\,\lambda_{\xi,e}\ne 0 \}$
 is finite and consists of nonzero elements pairwise different at $\pi_\xi$. Moreover,
 the representation is unique up to refinements of the partition of unity; cp.~\cite[\S\,6]{AK3} and
 \cite[\S\,5.1]{DOP}. The following result of~\cite[Proposition 4.6\,(1)]{Kus1}
 explains why and how the concept of local Hamel basis is a so useful technical tool (cp.~\cite{AK3}).

 \subsubsec{2.2.4.}~\proclaim{}Assume that  $\mathcal{E},\mathcal{X}\in\mathbb{V}^{(\mathbb B)}$,
 $[\![\,\mathcal{E}\subset\mathcal{X}\,]\!]\!=\!\mathbb{1}$,
 $[\![\,\mathcal{X}$ is a~vector subspace of
 $\mathcal{R}_{\mathbb{R}}]\!]=\mathbb{1}$, and\/ $X\!:=\!\mathcal{X}{\downarrow}$.
 Then $[\![\,\mathcal{E}$ is a~Hamel basis for the vector space $\mathcal{X}$
 $($over~$R^{\scriptscriptstyle\wedge})]\!]=\mathbb{1}$ if and only if $\mathcal{E}{\downarrow}$ is a~local Hamel basis for~$X$.
 \Endproc

 \subsubsec{2.2.5.}~A vector lattice $X$ is said to be {\it locally one-dimensional\/} if for every two nondisjoint $x_1,x_2\in X$ there exist nonzero components $u_1$ and $u_2$ of $x_1$ and $x_2$ respectively such that $u_1$ and $u_2$ are proportional; cp.~\cite[Definition~11.1]{AK3}. Equivalent definitions
see in~\cite[Proposition 5.1.2]{DOP}.

 \subsubsec{2.2.6.}~\proclaim{}Let $X$ be a laterally complete vector lattice and
 $\mathcal{X}\in\mathbb{V}^{(\mathbb{B})}$ be its Boolean valued representation
 with $\mathbb{B}\!:=\mathbb{P}(X)$. Them $X$ is locally one-dimensional if and only
 if $\mathcal{X}$ is one-dimensional vector lattice over $\mathbb{R}^{\scriptscriptstyle\wedge}$
 within $\mathbb{V}^{(\mathbb{B})}$, i.~e., $[\![\mathcal{R}=\mathbb{R}^{\scriptscriptstyle\wedge}]\!]=\mathbb{1}$.
  \Endproc

 \subsubsec{2.2.7.}~\proclaim{}A universally complete vector lattice is locally
 one-dimensional if and only if every band preserving linear operator in it is order bounded.
 \Endproc

 \beginproof~By the Gordon Theorem we can assume that $X=\mathcal{R}{\downarrow}$ with
 $\mathcal{R}\in\mathbb{V}^{(\mathbb{B})}$ and $\mathbb{B}\simeq\mathbb{P}(X)$. Thus, the
 problem reduces to existence of a~discontinuous solution to the Cauchy functional equation
 $(L)$ and the claim follows from 2.1.6.~\endproof

 \subsubsec{2.2.8.}~\proclaim{}Let $\mathbb{R}$ is a transcendental extension of a subfield
 $\mathbb{P}\subset\mathbb{R}$. There exists an $\mathbb{P}$-linear subspace $\mathcal{X}$
 in $\mathbb{R}$ such that $\mathcal{X}$ and $\mathbb{R}$ are isomorphic vector spaces
 over~$\mathbb{P}$ but they are not isomorphic as ordered vector spaces over $\mathbb{P}$.
 \Endproc

 \beginproof~Let $\mathcal{E}$ be a Hamel basis of a $\mathbb{P}$-vector space $\mathbb{R}$.
 Since $\mathcal{E}$ is infinite, we can choose a proper subset $\mathcal{E}_0\varsubsetneq\mathcal{E}$
 of the same cardinality: $|\mathcal{E}_0|=|\mathcal{E}|$. If $\mathcal{X}$ denotes the
 $\mathbb{P}$-subspace of $\mathbb{R}$ generated by $\mathcal{E}_0$,
 then $\mathcal{X}_0\varsubsetneq\mathbb{R}$ and $\mathcal{X}$ and $\mathbb{R}$ are isomorphic
 as vector spaces over $\mathbb{P}$. If $\mathcal{X}$ and $\mathbb{R}$ were isomorphic as ordered
 vector spaces over $\mathbb{P}$, then $\mathcal{X}$ would be order complete and, in consequence,
 we would have $\mathcal{X}=\mathbb{R}$, a contradiction.~\endproof

 \subsubsec{2.2.9.}~\proclaim{Theorem.}Let $X$ be a nonlocally one-dimensional universally complete
 vector lattice. Then there exist a vector sublattice $X_0\subset X$ and a band preserving linear bijection
 $T:X_0\to X$ such that $T^{-1}$ is also band preserving but $X_0$ and $X$ are not lattice isomorphic.
 \Endproc

 \beginproof~We can assume without loss of generality that $X=\mathcal{R}{\downarrow}$ and
 $[\![\mathcal{R}\ne \mathbb{R}^{\scriptscriptstyle\wedge}]\!] =\mathbb{1}$. By~4.6.5 there
 exist an $\mathbb{R}^{\scriptscriptstyle\wedge}$-linear subspace $\mathcal{X}$ in $\mathcal{R}$
 and $\mathbb{R}^{\scriptscriptstyle\wedge}$-linear isomorphism $\tau$
 from $\mathcal{X}$ onto $\mathcal{R}$, while $\mathcal{X}$ and $\mathcal{R}$ are not isomorphic
 as ordered vector spaces over $\mathbb{R}^{\scriptscriptstyle\wedge}$. Put
 $X_0\!:=\mathcal{X}{\downarrow}$, $T\!:=\tau{\downarrow}$ and $S\!:=\tau^{-1}{\downarrow}$.
 The maps $S$ and $T$ and are band preserving and linear. Moreover, $S=(\tau{\downarrow})^{-1}=T^{-1}$.
 It remains to observe that $X_0$ and $X$ are lattice isomorphic if and only if $\mathcal{X}$ and
 $\mathcal{R}$ are isomorphic as ordered vector spaces.~\endproof



 \subsubsec{2.2.10.}~Let $\gamma$ be a cardinal. A vector lattice $X$ is said to be
 {\it Hamel $\gamma$-ho\-mo\-ge\-neous\/} whenever there exists a local Hamel basis of
 cardinality $\gamma$ in $X$ consisting of strongly distinct weak order units. Two elements
 $x,y\in X$ are said to be {\it strongly distinct\/} if $|x-y|$ is a weak order unit in $X$.

 \subsubsec{2.2.11.}~\proclaim{}Let $X$ be a universally complete vector lattice. There is a band
 $X_0$ in $X$ such that $X_0^\perp$ is locally one-dimensional and there exists a partition of unity
 $(\pi_\gamma)_{\gamma\in\Gamma}$ in $\mathbb{P}(X_0)$ with $\Gamma$ a~set of infinite cardinals such
 that $\pi_\gamma X_0$ is Hamel $\gamma$-homogeneous for all $\gamma\in\Gamma$.
 \Endproc

 \subsubsec{2.2.12.}~A local Hamel basis is also called a $d$-\textit{basis}. This concept stems
 from~\cite{Coop}, but for the first time in the context of disjointness preserving operators was
 used in \cite{AVK, AVK1}. Various aspects of the concept can be found in \cite{AK3, AK2}.
 Theorem~2.2.6 was established in \cite{Gut6}, while 2.2.7 in \cite{AVK1, MW}.
 Another proof of Theorem~2.2.8 one can find in \cite{AK2}.

 \section{2.3. Algebraic Band Preserving Operators}

 In this section some description of algebraic orthomorphisms on a vector lattice is given and
 the Wickstead problem for algebraic operators is examined.

 \subsubsec{2.3.1.}~Let $\mathbb{P}[x]$ be a ring of polynomials in variable $x$ over a
 field $\mathbb{P}$. An operator $T$ on a vector space $X$ over a field $\mathbb{P}$ is said to
 be {\it algebraic\/} if there exists a~nonzero $\varphi\in\mathbb{P}[x]$, a polynomials with
 coefficients in $\mathbb{P}$, for which $\varphi(T)=0$.

 For an algebraic operator $T$ there exists a unique polynomial $\varphi_T$ such that $\varphi_T(T)=0$,
 the leading coefficient of $\varphi_T$ equals to 1, and $\varphi_T$ divides each polynomial $\psi$ with
 $\psi(T)=0$. The polynomial $\varphi_T$ is called the {\it minimal polynomial} of $T$. The simple
 examples of algebraic operators yield a projection $P$ (an idempotent operator, $P^2=P$) in $X$
 with $\varphi_P(\lambda)=\lambda^2-\lambda$ whenever $P\ne0,I_X$, and a nilpotent operator $S$
 ($S^m=0$ for some $m\in\mathbb{N}$) in $X$ with $\varphi_S(\lambda)=\lambda^k$, $k\leq m$.

 For an operator $T$ on $X$, the set of all eigenvalues of $T$ will be denoted throughout by $\sigma_p(T)$.
 A real number $\lambda$ is a root of $\varphi_T$ if and only if $\lambda\in\sigma_p(T)$.
 In particular, $\sigma_p(T)$ is finite. If $b-a^2>0$ for some $a,b\in\mathbb{R}$ then $T^2 + 2aT + bI$
 is a weak order unit in $\Orth(X)$ for every $T\in\Orth(X)$; cp.~\cite{BBS}.

 \subsubsec{2.3.2.}~\proclaim{}Let $X$ be a vector lattice and $T$ in $\Orth(X)$ is algebraic. Then
 $$
 \varphi_T(x)=\prod_{\lambda\in\sigma_p(T)}(x-\lambda).
 $$
 \Endproc

 \subsubsec{2.3.3.}~\proclaim{}Consider the universally complete vector lattice $X=\mathcal{R}{\downarrow}$. Let $T$
 be a~band preserving linear operator on $X$ and $\tau$ an $\mathbb{R}^{\scriptscriptstyle\wedge}$-linear
 function on $\mathcal{R}$. For $\varphi\in\mathbb{R}[x]$, $\varphi(x)=a_0+a_1x+\dots+a_nx^n$ define
 $\hat{\varphi}\in\mathbb{R}^{\scriptscriptstyle\wedge}[x]$ by
 $\hat{\varphi}(x)=a_0^{\scriptscriptstyle\wedge}+a_1^{\scriptscriptstyle\wedge}x
 +\dots+a_n^{\scriptscriptstyle\wedge}x^{n^{\scriptscriptstyle\wedge}}$. Then
 $$
 \hat{\varphi}(\tau){\downarrow}=\varphi(\tau{\downarrow}),\quad
 \varphi(T){\uparrow}=\hat{\varphi}(T{\uparrow}).
 $$
 \Endproc

 \beginproof~It follows from 1.2.4 and 1.2.7 that $(\tau^{n^{\scriptscriptstyle\wedge}}){\downarrow}=(\tau{\downarrow})^n$
 and $(T^{n}){\uparrow}=(T{\uparrow})^{n^{\scriptscriptstyle\wedge}}$. Thus, it remains to apply 1.4.9.~\endproof

 \subsubsec{2.3.4.}~A linear operator $T$ on a vector lattice $X$ is said to be {\it diagonal}
 if $T=\lambda_1 P_1+\dots+\lambda_m P_m$ for some collections of reals $\lambda_1,\dots,\lambda_m$
 and projection operators  $P_1,\dots,P_m$ on $X$ with $P_i\circ P_j=0$ $(i\ne j)$.
 In the equality above, we may and will assume that $P_1+\dots+P_n =I_X$ and that $\lambda_1,\dots,\lambda_m$
 are pairwise different. An algebraic operator $T$ is diagonal if and only if the minimal polynomial of $T$
 have the form $\varphi_T(x)=(x-\lambda_1)\cdot\ldots\cdot(x-\lambda_m)$ with pairwise different $\lambda_1,\dots,\lambda_m\in\mathbb{R}$.

 We call an operator $T$ on $X$  {\it strongly diagonal\/} if there exist pairwise disjoint band
 projections $P_1,\dots,P_m$ and real numbers $\lambda_1,\dots,\lambda_m$ such that $T=\lambda_1P_1+\dots+\lambda_mP_m$.
 In particular, every strongly diagonal operator on $X$ is an orthomorphism.

 \subsubsec{2.3.5.}~\proclaim{}Let $T=\lambda_1P_1+\dots+\lambda_mP_m$ be a diagonal operator on a vector lattice~$X$.
 Then $T$ is band preserving if and only if the projection operators $P_1,\dots,P_m$ are band preserving.
 \Endproc

 \beginproof~The sufficiency is obvious. To prove the necessity, observe first that if $T$ is band preserving then so is $T^n$  for all $n\in\mathbb{N}$ and thus $\varphi(T)$ is band preserving for every polynomial $\varphi\in\mathbb{R}[x]$.  Next, make use of the representation $P_j=\varphi_j(T)$ $(j\!:=1,\dots,m)$, where $\varphi_j\in\mathbb{R}[x]$ is an interpolation polynomial defined by $\varphi_j(\lambda_k)=\delta_{jk}$ with $\delta_{jk}$  the Kronecker symbol.~\endproof

 \subsubsec{2.3.6.}~\proclaim{Theorem.}Let $X$ be a universally complete vector lattice. The following
  are equivalent:

 \subsubsec{(1)}~The Boolean algebra $\mathbb{P}(X)$ is
 $\sigma$-distributive.

 \subsubsec{(2)}~Every algebraic band preserving operator in $X$ is
 order bounded.

 \subsubsec{(3)}~Every algebraic band preserving operator in $X$ is
 strongly diagonal.


 \subsubsec{(4)}~Every band preserving diagonal operator in $X$ is
 strongly diagonal.


 \subsubsec{(5)}~Every band preserving nilpotent operator in $X$ is
 order bounded.

 \subsubsec{(6)}~Every band preserving nilpotent operator in $X$ is trivial.
 \Endproc

 \beginproof~The only nontrivial implications are (2)~$\Longrightarrow$~(3) and (6)~$\Longrightarrow$~(2).

 (2)~$\Longrightarrow$~(3)~We have to prove that an algebraic orthomorphism on $X$ is strongly diagonal. Let $T$ be an orthomorphism in $X$ and $\varphi(T)=0$, where $\varphi$ is a minimal polynomial of $T$, so that $\varphi(\lambda)=(\lambda-\lambda_1)\cdot\ldots\cdot(\lambda-\lambda_m)$ with $\lambda_1,\dots,\lambda_m\in\mathbb{R}$.
 Since $T$ admits a unique extension to an orthomorphism on $X^{\rm u}$, we can assume without loss of generality that $X=X^{\rm u}=\mathcal{R}{\downarrow}$ and $\tau=T{\uparrow}$. Then $[\![\tau(x)=\lambda_0 x\ (x\in\mathcal{R})]\!]=\mathbb{1}$ for some $\lambda_0\in\mathcal{R}$. It is seen from  2.3.3 that $\hat{\varphi}(\lambda_0)=0$ and thus $(\lambda_0-\lambda_1^{\scriptscriptstyle\wedge})\cdot\ldots\cdot(\lambda_0-\lambda_m^{\scriptscriptstyle\wedge})=0$
 or $\lambda_0\in\{\lambda_1^{\scriptscriptstyle\wedge},\dots,\lambda_m^{\scriptscriptstyle\wedge}\}$ within $\mathbb{V}^{(\mathbb{B})}$.
 Put $P_l\!:=\chi(b_l)$ with $b_l\!:=[\![\lambda_0=\lambda_l^{\scriptscriptstyle\wedge}]\!]$
 and observe that $\{P_1,\dots,P_m\}$ is a partition of unity in $\mathbb{P}(X)$. Moreover, given $x\in X$, we see that $b_l\leq[\![Tx=\tau x=\lambda_0 x]\!]
 \wedge[\![\lambda_0=\lambda_l^{\scriptscriptstyle\wedge}]\!]\leq[\![Tx=\lambda_l^{\scriptscriptstyle\wedge}x]\!]$,
 so that $P_l Tx=P_l(\lambda_l x)=\lambda_lP_l(x)$. Summing up over $l=1,\dots,m$,
 we get $Tx=\lambda_1P_1x+\dots+\lambda_mP_m$.

 (6)~$\Longrightarrow$~(1)~Arguing for a contradiction, assume that 2.3.6\,(2) is fulfilled and construct a nonzero band preserving nilpotent operator in $X$. By~2.2.7 and 2.1.6 we have $\mathbb{V}^{(\mathbb{B})}\models\mathcal{R}\ne\mathbb{R}^{\scriptscriptstyle\wedge}$
 and thus $\mathcal{R}$ is an~infinite-dimensional
 vector space over~$\mathbb{R}^{\scriptscriptstyle\wedge}$ within $\mathbb{V}^{(\mathbb{B})}$. Let $\mathcal{E}\subset\mathcal{R}$ be a Hamel basis and choose
 an infinite sequence $(e_n)_{n\in\mathbb{N}}$ of pairwise distinct elements in $\mathcal{E}$. Fix a natural $m>1$ and define an $\mathbb{R}^{\scriptscriptstyle\wedge}$-linear
 function $\tau:\mathcal{R}\to\mathcal{R}$ within $\mathbb{V}^{(\mathbb{B})}$
 by letting $\tau(e_{km+i})=e_{km+i-1}$ if $2\leq i\leq m$,
 $\tau(e_{km+1})=0$ for all $k\!:=0,1,\dots$, and $\tau(e)=0$ if $e\ne e_n$
 for all $n\in\mathbb{N}$. In other words, if $\mathcal{R}_0$ is the
 $\mathbb{R}^{\scriptscriptstyle\wedge}$-linear subspace of
 $\mathcal{R}$ generated by the sequence $(e_n)_{n\in\mathbb{N}}$,
 then $\mathcal{R}_0$ is an invariant subspace for $\tau$ and
 $\tau$ is the linear operator associated to the infinite
 block matrix $\diag(A,\dots,A,\dots)$ with equal blocks in the
 principal diagonal and $A$ a~the Jordan block of size $m$ with eigenvalue $0$.
 It follows that $\tau$ is discontinuous and $\tau^m=0$ by construction. Consequently,
 $T\!:=\tau{\downarrow}$ is a band preserving linear operator in $X$ and $T^m=0$ by  2.3.3,
 but $T$ is not order bounded; a contradiction.~\endproof

 \subsubsec{2.3.7.}~Algebraic order bounded disjointness preserving operators in vector lattices were treated in~\cite{BBS} where, in particular, the Propositions~2.3.2 and 2.3.5 were proved. Theorem 2.3.6 was obtained in \cite{KZ}.

 \section{2.4. Involutions and Complex Structures}

 The main result of this section tells us that in a real non locally one-dimensional universally
 complete vector lattice there are band preserving complex structures and nontrivial band preserving involutions.

 \subsubsec{2.4.1.}~A linear operator $T$ on a vector lattice $X$ is called {\it involutory} or
 an {\it involution} if\/ $T\circ T=I_X$ (or, equivalently, $T^{-1}=T$) and is called a {\it complex
 structure\/} if $T\circ T=-I_X$ (or, equivalently, $T^{-1}=-T$). The operator $P-P^\perp$, where $P$ is
 a~projection operator on $X$ and $P^\perp=I_X-P$, is an involution. The involution $P-P^\perp$ with band
 projections $P$ is referred to as {\it trivial\/}.

\subsubsec{2.4.2.}\proclaim{}~Let $X$ be a Dedekind complete vector
lattice. Then there is no order bounded band preserving complex
structure in $X$ and there is no nontrivial order bounded band
preserving involution in $X$. \Endproc

\beginproof~An order bounded band preserving operator $T$ on a universally complete vector
lattice $X$ with weak unit $\mathbb{1}$ is a multiplication
operator: $Tx=ax$ $(x\in X)$ for some $a\in X$. It follows that $T$
is an involution if and only if  $a^{2}=\mathbb{1}$ and hence there
is a band projection $P$ on $E$ with
$a=P\mathbb{1}-P^\perp\mathbb{1}$ or $T=P-P^{\bot}$. If $T$ is a
complex structure on $E$ then the corresponding equation
$a^2=-\mathbb{1}$ has no solution.~\endproof

 \subsubsec{2.4.3.}~\proclaim{Theorem.}~Let $\mathbb{F}$ be a dense subfield of\/
 $\mathbb{R}$ and $B\subset \mathbb{R}$ be a non-empty finite or countable set. Then there exists a
 discontinuous $\mathbb{F}$-linear function $f:\mathbb{R}\rightarrow\mathbb{R}$ such
 that $f\circ f=f$ and $f(x)=x$ for all $x\in B$.
 \Endproc

 \beginproof~Let $\mathcal{E}\subset \mathbb{R}$ be a Hamel basis of $\mathcal{R}$ over
 $\mathbb{R}^{\scriptscriptstyle\wedge}$. Every $x\in B$ can be written in the form
 $x=\sum_{e\in\mathcal{E}}\lambda_{e}(x)e$, where $\lambda_e(x)\in\mathbb{F}$ for all
 $e\in\mathcal{E}$. Put ${\mathcal{E}(x)\!:=\{e\in\mathcal{E}:\,\lambda_e(x)\ne0\}}$ and
 $\mathcal{E}_{0}=\bigcup_{x\in B}\mathcal{E}(x)$. Since $B$ is finite or countable, so is also $\mathcal{E}_{0}$,
 Hence $\mathcal{E}\backslash \mathcal{E}_{0}$ has the cardinality of continuum. There exists a decomposition
 $\mathcal{E}_{1}\cup\mathcal{E}_{2}=\varnothing$,
 where $\mathcal{E}_{1}$ and $\mathcal{E}_{2}$ disjoint sets
 both having the same cardinality. Hence there exists a one-to-one mapping $g_{0}$ from $\mathcal{E}_{1}$ onto $\mathcal{E}_{2}$ with the inverse $g_{0}^{-1}:\mathcal{E}_{2}\to\mathcal{E}_{1}$.

 Now we define the function $g:\mathcal{E}\rightarrow\mathcal{E}$ as follows:
 \begin{equation}\label {G}
 g(h)=
 \begin{cases}
 g_{0}(h),&\text{~~for}~h\in\mathcal{E}_{1},
 \\
 g^{-1}_{0}(h),&\text{~~for}~h\in\mathcal{E}_{2},
 \\
 h,&\text{~~for~}~h\in\mathcal{E}_{0}.
 \end{cases}
 \end{equation}
 It can be easily checked that the $\mathbb{F}$-linear extension $f:\mathbb{R}\rightarrow\mathbb{R}$
 od a function $g$ is the sought involution.~\endproof

 \subsubsec{2.4.4.}~\proclaim{Theorem.}~Let $\mathbb{F}$ be a dense subfield of\/ $\mathbb{R}$.
 Then there exists a~dis\-con\-tinuous $\mathbb{F}$-linear function $f:\mathbb{R}\rightarrow\mathbb{R}$
 such that $f\circ f=-f$.
 \Endproc

 \beginproof~The proof is similar to that of Theorem 4.13.3
 with the minor modifications: put $\mathcal{E}_0=\varnothing$ and define
 \begin{equation*}
 g(h)=
 \begin{cases}
 -g_{0}(h),&\text{~~for}~h\in\mathcal{E}_{1},
 \\
 g^{-1}_{0}(h),&\text{~~for}~h\in\mathcal{E}_{2}.
 \end{cases}\ \ \endproof
 \end{equation*}

 Interpreting Theorems 2.4.3 and 2.4.4 in a Boolean valued model yields the   result.

 \subsubsec{2.4.5.}~\proclaim{Theorem.}~Let $X$ be a universally complete real vector lattice that is not locally one-dimensional. Then
 \vspace{2pt}

 \subsubsec{(1)}~For every nonempty finite or countable set $B\subset X$ there exists a band preserving involution $T$ on $X$ with $T(x)=x$ for all $x\in B$.%
 \vspace{2pt}

 \subsubsec{(2)}~There exists a band preserving complex structure on $X$.
 \Endproc

 \beginproof~Assume that $X=\mathcal{R}{\downarrow}$. Take a one-to-one function
 $\nu:\mathbb{N}\to X$ with $B=\im(\nu)$. The function $\nu{\upwardarrow}:\mathbb{N}^
 {\scriptscriptstyle\wedge}\to X$ may fail to be one-to-one within $\mathbb{V}^{(\mathbb{B})}$
 but $B{\uparrow}$ is again finite or countable, as  $B{\uparrow}=\im(\nu{\upwardarrow})$
 by 1.2.7. By Theorem 2.4.3 there exists an $\mathbb{R}^{\scriptscriptstyle\wedge}$-linear
 function $\tau:\mathcal{R}\to\mathcal{R}$
 such that $[\![\tau\circ\tau=I_\mathcal{R}]\!]=\mathbb{1}$ and
 $$
 \gathered
 \mathbb{1}=[\![(\forall\, x\in B{\uparrow})\tau(x)=x]\!]=
 [\![(\forall\,n\in\mathbb{N}^{\scriptscriptstyle\wedge})
 \tau(\nu{\upwardarrow}(n)=\nu{\upwardarrow}(n)]\!]
 \\
 =\bigwedge_{n\in\mathbb{N}}[\![\tau(\nu{\upwardarrow}(n^{\scriptscriptstyle\wedge}))
 =\nu{\upwardarrow}(n^{\scriptscriptstyle\wedge})]\!]=
 \bigwedge_{n\in\mathbb{N}}[\![\tau(\nu(n))=\nu(n)]\!]
 \\
 =\bigwedge_{n\in\mathbb{N}}[\![\tau{\downarrow}(\nu(n))=
 \nu(n)]\!].
 \endgathered
 $$
 It follows that if $T\!:=\tau{\downarrow}$ then $T\circ T=I_X$ by 1.2.4 and  $T(\nu(n))=\nu(n)$ for all $n\in\mathbb{N}$ as required in (1). The second claim is proved in a similar way on using Theorem 2.4.4.~\endproof

 \subsubsec{2.4.6.}~\proclaim{Corollary.}~Let $X$ be a universally complete vector lattice. Then the following are equivalent:
 \vspace{2pt}

 \subsubsec{(1)}~$X$ is locally one-dimensional.
 \vspace{2pt}

 \subsubsec{(2)}~There is no nontrivial band preserving involution on $X$.

 \vspace{2pt}
 \subsubsec{(3)}~There is no band preserving complex structure on $X$.
 \Endproc

 \subsubsec{2.4.7.}~\proclaim{Corollary.}~Let $X$ be a universally complete real vector lattice.
 Then $X$ admits a structure of complex vector space with a band preserving complex multiplication.
 \Endproc

 \beginproof~A complex structure $T$ on $X$ allows to define on $X$ a structure of a vector space over
 the complexes $\mathbb{C}$, by setting $(\alpha + i\beta)x =
 \alpha x+\beta T(x)$ for all $z=\alpha+i\beta\in\mathbb{C}$ and $x\in E$. If $T$ is band preserving
 then the map $x\mapsto zx$ $(x\in E)$ is evidently band preserving for all $z\in\mathbb{C}$.~\endproof

 \subsubsec{2.4.8.} The main results of this section were obtained in \cite{KZ1}. In connection with Corollary 2.4.7  spaces without complex structure should be mentioned, see  \cite{Die, GM1, GM2, Sza}.

 \section{Chapter 3. Disjointness Preserving Operators}

 \subsection{3.1. Characterization and Representation}

 Now we will demonstrate that some properties of disjoint preserving operators
 are just Boolean valued interpretations of elementary properties of~dis\-jointness
 preserving functionals.

 \subsubsec{3.1.1.}~\theorem{}Assume that $Y$ has the projection property. An order bounded
 linear operator $T:X\to Y$ is disjointness preserving if and only if $\ker(bT)$ is an order
 ideal in $X$ for every projection $b\in\mathbb{P}(Y)$.
 \Endproc

 \beginproof~The necessity is obvious, and so only the sufficiency will be proved. Suppose that
 $\ker(bT)$ is an order ideal in $X$ for every $b\in\mathbb{P}(Y)$. We can assume that
 $Y\subset\mathcal{R}{\downarrow}$ by the Gordon Theorem. Let $|y|\leq|x|$ and $b\!:=[\![Tx=0]\!]$.
 Then $bTx=0$ by $(G)$ and $bTy=0$ by the hypothesis. Again, making use of $(G)$ we have $b\leq[\![Ty=0]\!]$.
 Thus $[\![Tx=0]\!]\leq[\![Ty=0]\!]$ or, what is the same, $[\![Tx=0]\!]\Rightarrow[\![Ty=0]\!]=\mathbb{1}$.
 Now, put $\tau\!:=T{\upwardarrow}$ and ensure that
 $\ker(\tau)$ is an order ideal in $X^{\scriptscriptstyle\wedge}$
 within~$\mathrm{V}^{(\mathbb{B})}$. Making use of the fact that
 $|x|\leq|y|$ if and only if $[\![x^{\scriptscriptstyle\wedge}\leq
 y^{\scriptscriptstyle\wedge}]\!]=\mathbb{1}$, we deduce:
 $$
 \gathered
 {}[\![\ker(\tau)\text{~is an order ideal in~}
 X^{\scriptscriptstyle\wedge}]\!]
 \\
 =[\![(\forall\,x,y\in X^{\scriptscriptstyle\wedge})\,
 (\tau(x)=0\wedge|y|\leq|x|\,\rightarrow\,\tau(y)=0)]\!]
 \\
 =\bigwedge_{x,y\in X}[\![(\tau(x^{\scriptscriptstyle\wedge})=0]\!]
 \wedge[\![|y^{\scriptscriptstyle\wedge}|\leq|x^{\scriptscriptstyle\wedge}|]\!]\,
 \Rightarrow\,[\![\tau(y^{\scriptscriptstyle\wedge})=0]\!]
 \\
 =\bigwedge\Big\{[\![T(x)=0]\!]\,\Rightarrow\,[\![T(y)=0]\!]:\ x,y\in X,\
 |y|\leq|x|\Big\}=\mathbb{1}.
 \endgathered
 $$
 Apply within~$\mathbb{V}^{(\mathbb{B})}$ the fact that the functional $\tau$ is disjointness preserving if and only if $\ker(\tau)$ is an order ideal in $X^{\scriptscriptstyle\wedge}$. It follows that $T$ is also disjointness preserving according to 1.4.4\,(5).~\endproof

 \subsubsec{3.1.2.}~Similar reasoning yields that if $Y$ has the projection property then for an~order bounded disjointness preserving linear operator $T\in L^{\scriptscriptstyle\sim} (X,Y)$ there exists a~band projection $\pi\in\mathbb{P}(Y)$
 such that $T^{+}=\pi|T|$ and $T^{-}=\pi^\perp|T|$. In particular, $T=(\pi-\pi^\perp)|T|$ and $|T|=(\pi-\pi^\perp)T$. To ensure this, observe that the functional $\tau\!:=T{\upwardarrow}$ preserves disjointness if and only if either $\tau$, or $-\tau$ is a lattice homomorphism.

 From this fact it follows that $T\in L^{\scriptscriptstyle\sim}(X,Y)$ is disjointness preserving if and only if $(Tx)^+\perp (Ty)^-$ for all $x,y\in X_+$. Indeed, for arbitrary $x,y\in X_+$ we can write $(Tx)^+=(Tx)\vee0\leq T^+x=\pi|T|x$ and, similarly,
 $(Ty)^-\leq \pi^\perp|T|y$. Hence $(Tx)^+\wedge (Ty)^-=0$.

 \subsubsec{3.1.3.}~\Theorem{}Let $X$ and $Y$ be vector lattices with $Y$ Dedekind
 complete. For a pair of disjointness preserving operators $T_1$ and $T_2$ from
 $X$ to $Y$ there exist a band projection $\pi\in\mathbb{P}(Y)$, a lattice
 homomorphism $T\in\Hom(X,Y)$, and ortho\-mor\-phisms $S_1,S_2\in\Orth(Y)$ such that
 $$
 \gathered
 |S_1|+|S_2|=\pi,\ \ \pi T_1=S_1T,\ \ \pi T_2=S_2T ,
 \\
 \im(\pi^\perp T_1)^{\perp\perp}=
 \im(\pi^\perp T_2)^{\perp\perp}=\pi(Y),\ \
 \pi^\perp T_1\perp\pi^\perp T_2.
 \endgathered
 $$
 \Endproc

 \beginproof~As usual, there is no loss of generality in assuming that ${Y=\mathcal{R}{\downarrow}}$. Put $\tau_1\!:=T_1{\upwardarrow}$ and $\tau_2\!:=T_2{\upwardarrow}$. The desired result is a Boolean valued interpretation of the following fact: If the disjointness preserving functionals $\tau_1$ and $\tau_2$ are not proportional then they are nonzero and disjoint. Put
 $b\!:=[\![\tau_1\text{~and~}\tau_2$ are proportional$]\!]$ and $\pi\!:=\chi(b)$. Then within $\mathbb{V}^{([\mathbb{0},b])}$
 there exist a lattice homomorphism $\tau:X^{\scriptscriptstyle\wedge}\to\mathcal{R}$ and reals $\sigma_1,\sigma_2\in\mathcal{R}$ such that $\tau_i=\sigma_i\tau$. If the function $\bar{\sigma}_i$ is defined as $\bar{\sigma}_i(\lambda)\!:=\sigma_i\lambda$ $(\lambda\in\mathcal{R})$, then the operators  $S_1\!:=\sigma_1{\downwardarrow}$, $S_2\!:=\sigma_2{\downwardarrow}$, and $T\!:=\tau{\downwardarrow}$ satisfy the first line of required conditions. Moreover, $b^\ast=[\![\tau_1\ne0]\!]\wedge[\![\tau_1\ne0]\!]\wedge
 [\![|\tau_1|\wedge|\tau_2|=0]\!]$ and $\pi^\perp=\chi(b^\ast)$, so that the second line of required conditions is also satisfied.~\endproof

 \subsubsec{3.1.4.}~\Corollary{}Let $X$ and $Y$ be vector lattices with $Y$ Dedekind
 complete. The sum $T_1+T_2$ of two disjointness preserving operators $T_1,T_2:X\to Y$ is disjointness
 preserving if and only if there exist pairwise disjoint band projections $\pi,\pi_1,\pi_2\in\mathbb{P}(Y)$,
 orthomorphisms $S_1,S_2\in\Orth(Y)$ and a lattice homomorphism $T\in\Hom(X,Y)$ such
 that the following system of relations is consistent
 $$
 \gathered
 \pi+\pi_1+\pi_2=I_Y,\ \ |S_1|+|S_2|=\pi,
 \\
 T(X)^{\perp\perp}=\pi(Y),\ \ \pi_1T_2=\pi_2T_1=0,
 \\
 \pi T_1=S_1T,\ \ \pi T_2=S_2T.
 \endgathered
 $$
 Consequently, in this case $T_1+T_2=\pi_1T_1+\pi_2T_2+(S_1+S_2)T$.
 \Endproc

 \subsubsec{3.1.5.}~\corollary{}The sum $T_1+T_2$ of two disjointness preserving operators
 $T_1,T_2:X\to Y$ is disjointness preserving if and only if $T_1(x_1)\perp T_2(x_2)$
 for all $x_1,x_2\in X$ with $x_1\perp x_2$.
 \Endproc

 \beginproof~The necessity is immediate from Theorem 3.1.4, since $T_1=\pi_1T_1+S_1T$
 and $T_2=\pi_2T_2+S_2T$. To see the sufficiency, observe that if $T_1$ and $T_2$ meet
 the above condition then $T_kx_1\perp T_lx_2$ $(k,l\!:=1,2)$ and thus $(T_1+T_2)(x_1)\perp(T_1+T_2)(x_2)$
 for any pair of disjoint elements $x_1, x_2\in X$.~\endproof

 \subsubsec{3.1.6.}~Aspects of the theory of disjointness preserving operators are presented
 in \cite{Gut1, DOP, DOR}. The recent results on disjointness preserving operators are surveyed in~\cite{Boul}.
 In particular, the concept of disjointness preserving set of operators is discussed in the survey.
 Using this concept Corollary 3.1.5 may be reformulated as follows:
 \textsl{$T_1+T_2$ is disjointness preserving if and only if $\{T_1,T_2\}$ is a disjointness preserving set of operators.}

 \subsection{3.2. Polydisjoint Operators}

 The~aim of the~present section is to describe the~order ideal
 in the~space of order bounded operators generated by order
 bounded disjointness preserving operators ($=$ $d$-homomorphisms)
 in terms of $n$-disjoint operators.

 \subsubsec{3.2.1.}~Let $X$ and $Y$ be vector lattices and $n$
 be a positive integer. A~linear operator~$T:X\rightarrow Y$
 is said to be {\it $n$-disjoint}  if, for any collection of
 $n+1$ pairwise disjoint elements $x_0,\ldots,x_n\in X$,
 the~infimum of the~set $\big\{|Tx_k|:\,k\!:=0,1,\ldots,n\big\}$
 equals zero; symbolically:
 $$
 (\forall\, x_0,x_1\ldots,x_n\in X)\ x_k\perp x_l\ (k\neq l)\
 \Longrightarrow\ |Tx_0|\wedge\ldots\wedge|Tx_n|=0.
 $$
 An operator is called \textit{polydisjoint} if it is $n$-disjoint for some $n\in\mathbb{N}$. A $1$-disjoint operator is just a disjointness preserving operator.

 \subsubsec{3.2.2.}~Consider some simple properties of $n$-disjoint operators. Let~$X$ and $Y$ be vector lattices
 with $Y$ Dedekind complete.

 \subsubsec{(1)}~An~operator $T\in L^{\scriptscriptstyle\sim}(X, Y)$ is $n$-disjoint if and only if $|T|$ is~$n$-disjoint.

 \subsubsec{(2)}~Let $T_1,\ldots,T_n$ be order bounded and disjointness preserving operators from $X$ to $Y$. Then the operator $T\!:=T_1+\ldots+T_n$ is $n$-disjoint.

 \subsubsec{3.2.3.}\proclaim{}An order bounded functional on a
 vector lattice is $n$-disjoint if and only if it is representable as a~disjoint sum of~$n$ order bounded disjointness preserving functionals. Such representation is unique up to  permutation.\Endproc

 \beginproof~Assume that $f$ is a positive $n$-disjoint functional on a vector lattice $C(Q)$. Prove that the corresponding Radon measure $\mu$ is the sum of $n$ Dirac measures. This is equivalent to saying that the
 support of $\mu$ consists of $n$ points. If~there are $n+1$
 points $q_0,q_1,\ldots,q_n\in Q$ in the support of $\mu$ then
 we may choose pair-wise disjoint compact neighborhoods
 $U_0,U_1,\ldots,U_n\in Q$ of these points and next take pair-wise disjoint open sets $V_k\subset Q$ with $\mu(U_k)>0$ and $U_k\subset V_k$ $(k=0,1,\ldots,n)$. Using the Tietze--Urysohn Theorem construct a continuous function $x_k$ on $Q$ which vanishes on $Q\setminus V_k$ and is identically equal to $1$ on $U_k$. Then $x_0\wedge x_1\wedge\ldots\wedge x_n=0$ but none of $f(x_0),f(x_1),\ldots,f(x_n)$ is equal to zero, since $f(x_k)\geq\mu(U_k)>0$ for all $k\!:=0,1,\ldots,n$. This contradiction shows that the support of $\mu$ consists of~$n$ points. The general case is reduced to what was proven by using the Kre\u\i ns--Kakutani Representation Theorem.~\endproof

 \subsubsec{3.2.4.}~\theorem{}An~order bounded operator from a
 vector lattice to a Dedekind complete vector lattice is
 $n$-disjoint for some $n\in\mathbb{N}$ if and only if~it is
 representable as a~disjoint sum of~$n$ order bounded
 disjointness preserving operators.\Endproc

 \beginproof~Assume that the operator $T\in L^\sim(X,Y)$ is  \mbox{$n$-disjiint} and denote $\tau\!:=T\upwardarrow\in\mathbb{V}^{(\mathbb{B})}$. It is deduced by direct calculation of Boolean truth values that  $\tau:X^{\scriptscriptstyle\wedge}\to\mathcal{R}$ is an order bounded \mbox{$n$-disjoint} functional within $\mathbb{V}^{(\mathbb{B})}$. According to Transfer Principle, applying Proposition 3.2.3 to $\tau$ yields pairwise disjoint order bounded disjointness preserving functionals $\tau_1,\ldots,\tau_n$ on
 $X^{\scriptscriptstyle\wedge}$ with $\tau=\tau_1+\ldots+\tau_n$.
 It remains to observe that the linear operators $T_1\!:=\tau_1{}\downwardarrow,\dots,$
 $T_n\!:=\tau_n{}\downwardarrow$ from $X$ to $Y$  are order bounded, disjointness preserving, and $T_1+\ldots+T_n=T$. Moreover, if $k\ne j$ then
 $0=(\tau_k\wedge\tau_l){\downwardarrow}=
 \tau_k{\downwardarrow}\wedge\tau_l{\downwardarrow}=T_k\wedge T_l$, so that $T_k$ and $T_l$ are disjoint.~\endproof

  \subsubsec{3.2.5.}~It can be easily seen that the representation of an order bounded  $n$-dis\-joint
  operator in Theorem $3.1.4$ is unique up to mixing: if  $T=T_1+\ldots+T_n=S_1+\ldots+S_m$ for two pairwise
  disjoint collections $\{T_1,\ldots,T_n\}$ and $\{S_1,\ldots,T_n\}$ of order bounded disjointness  preserving
  operators then for every $j\!=1,\ldots,m$ there exists a disjoint family of projections
  $\pi_{1j},\ldots,\pi_{nj}\in\mathbb{P}(Y)$ such that  $S_j=\pi_{1j}T_1+\ldots+\pi_{nj}T_n$ for all $j\!:=1,\ldots,m$.

 \subsubsec{3.2.6.}~\corollary{}A positive operator operator from a vector lattice to a Dedekind complete vector lattice is
 $n$-disjoint if and only if it is the sum of $n$ lattice
 homomorphisms. \Endproc

 \subsubsec{3.2.7.}~\corollary{}The set of polydisjoint operators from a vector lattices to a Dedekind
 complete vector lattices coincides with the order ideal in the vector lattice of order bounded operators
 generated by lattice homomorphisms. \Endproc

 \subsubsec{3.2.8.}~The characterizations of sums of lattice homomorphisms (Corollary~3.2.6), sums of
 disjointness preserving operators (Theorem 3.2.4), and the ideal of order bounded operators generated by lattice
 homomorphism (Corollary~3.2.7) were proved in \cite{BHP} using standard tools.
 An algebraic approach to the problem see in \cite{DOP}.

 \subsection{3.3. Differences of Lattice Homomorphisms}

 This section presents a characterization of order bounded operators representable as a difference of two lattice homomorphisms. The starting point of this question is the celebrated Stone Theorem about the structure of vector sublattices in the Banach lattice $C(Q,\mathbb R)$ of continuous real functions on a compact space $Q$. This theorem may
 be rephrased in the above terms as follows:

 \subsubsec{3.3.1.}~\proclaim{Stone Theorem.}Each closed vector sublattice of $C(Q,\mathbb R)$ is the intersection of the kernels of some differences of lattice homomorphisms on~$C(Q,\mathbb{R})$.
 \Endproc

 \subsubsec{3.3.2.}~In view of this theorem it is reasonable to
 refer to a difference of lattice homomorphisms on a~vector
 lattice~$X$ as a~two-point relation on~$X$. We are not obliged to assume here that the lattice homomorphisms under study act into the reals~$\mathbb R$. Thus a linear operator $T:X\to Y$ between vector lattices is said to be a~{\it two-point relation\/} on $X$ whenever it is written as a difference of two lattice homomorphisms. An operator $bT\!:=b\circ T$ with $b\in\mathbb{B}\!:=\mathbb{P}(Y)$ is called a~{\it stratum\/} of~$T$.

 \subsubsec{3.3.3.}~The kernel $\ker(bT)$ of any stratum of a two point relation $T$ is evidently a sublattice of $X$, since it is determined by an equation $bT_1x=bT_2x$. Thus, each stratum $bT$ of an~order bounded disjointness preserving operator $T:X\to Y$ is a two-point relation on $X$ and so its kernel is a vector sublattice of~$X$. The main result of this section says that the converse is valid too. To handle the corresponding scalar problem a formula of subdifferential calculus is used; cp.~ \cite{Kut2, KK3}. In the following form of this auxiliary fact \textit{positive decomposition\/} of a functional $f$ means any representation $f=f_1+\ldots+f_N$ with positive functionals $f_1,\ldots,f_N$.

 \subsubsec{3.3.4.}~\proclaim{Decomposition Theorem.}Assume that $H_1,\ldots,H_N$ are cones in a~vector lattice $X$
 and $f$ and $g$ are positive functionals on $X$. The inequality
 $$
 f(h_1\vee\ldots\vee h_N)\ge g(h_1\vee\ldots\vee h_N)
 $$
 holds for all $h_k\in H_k$ $(k:=1,\ldots,N)$ if and only if to each positive decomposition $(g_1,\ldots,g_N)$ of $g$ there is a positive decomposition $(f_1,\ldots,f_N)$ of $f$ such that
 $$
 f_k(h_k)\ge g_k(h_k)\quad (h_k\in H_k;\ k:=1,\ldots,N).
 $$ \Endproc

 \subsubsec{3.3.5.}\proclaim{}Let $\mathbb{F}$ be a dense subfield in $\mathbb{R}$ and $X$ a vector lattice
 over $\mathbb{F}$. An order bounded $\mathbb{F}$-linear functional from $X$ to $\mathbb{R}$ is a two-point
 relation if and only if its kernel is a $\mathbb{F}$-linear sublattice of the ambient vector lattice. \Endproc

 \beginproof~Let $l$ be an order bounded functional on a vector lattice $X$. Denote $f\!:=l^+$, $g\!:=l^-$,
 and $H\!:=\ker(l)$. It suffices to demonstrate only that~$g$ is a~lattice homomorphism, i.\,e., $[0,g]=[0,1]g$; cp.~ \cite{DOP}.
 So, we take $0\le g_1\le g$ and put $g_2:=g-g_1$. We~may assume that $g_1\ne0$ and $g_1\ne g$. By~hypothesis,
 for all $h_1,h_2\in\ker(l)$ we have the $f(h_1\vee h_2)\ge g(h_1\vee h_2)$. By the Decomposition
 theorem there is a positive decomposition $f=f_1+f_2$ such that $f_1(h)-g_1(h)=0$ and  $f_2(h)-g_2(h)=0$
 for all $h\in H$. Since $H=\ker(f-g)$, we see that there are reals $\alpha$ and $\beta$ satisfying
 $f_1-g_1=\alpha(f-g)$ and $f_2-g_2=\beta(f-g)$.
 Clearly,  $\alpha+\beta=1$ (for otherwise $f=g$ and $l=0$).
 Therefore, one of the reals $\alpha$ and $\beta$ is strictly
 positive. If $\alpha>0$ then we have $g_1=\alpha g$ for~$f$ and $g$ are disjoint. If $\beta>0$ then,
 arguing similarly, we see that $g_2=\beta g$. Hence, $0\le\beta\le1$ and we again see that $g_1\in[0,1]g$.~\endproof

 \subsubsec{3.3.4.}~\theorem{}An order bounded operator from a vector lattice to a~Dedekind
 complete vector lattice is a two-point relation if and only if the kernel of its every
 stratum is a vector sublattice of the ambient vector lattice. \Endproc

 \beginproof~The necessity is obvious, so only the sufficiency will be proved.
 Let $T\in L^\sim(X,Y)$ and $\ker(bT):=(bT)^{-1}(0)$ is a vector sublattice in $X$
 for all $b\in\mathbb{P}(Y)$. We apply the Boolean valued ``scalarization'' putting $Y=\mathcal{R}{\downarrow}$.

 Denote $\tau:=T{\upwardarrow}$ and observe that the validity of identities $T^+{\upwardarrow}=\tau^+$ and $T^-{\upwardarrow}=\tau^-$ within $\mathbb{V}^{(\mathbb{B})}$ is proved by easy calculation of Boolean truth values.
 Moreover, $[\![\ker(\tau)$ is a~vector sublattice of~${X^{\scriptscriptstyle{\wedge}} ]\!]=\mathbb{1}}$.
 Indeed, given $x,y\in X$, put
 $b:=[\![Tx=0^{\scriptscriptstyle{\wedge}}]\!]\wedge
 [\![Ty=0^{\scriptscriptstyle{\wedge}}]\!]$.
 This means that $x,y\in\ker(bT)$. Hence, we see by hypothesis that $bT(x\vee y)=0$, whence
 $b\leq[\![T(x\vee y)=0^{\scriptscriptstyle{\wedge}}]\!]$. Replacing $T$ by $\tau$ yields
 $[\![\tau(x^{\scriptscriptstyle{\wedge}})
 =0^{\scriptscriptstyle{\wedge}}\wedge
 \tau(y^{\scriptscriptstyle{\wedge}})=
 0^{\scriptscriptstyle{\wedge}}]\!]\leq
 [\![\tau(x\vee y)^{\scriptscriptstyle{\wedge}}=
 0^{\scriptscriptstyle{\wedge}}]\!]$.
 A straightforward calculation of Boolean truth values completes the proof:
 $$
 [\![\ker(\tau)\ \text{\rm is a\ Riesz\ subspace\ of}\
 X^{\scriptscriptstyle{\wedge}}]\!]
 $$
 $$
 \gathered
 =[\![(\forall\,x,y\in X^{\scriptscriptstyle{\wedge}})
 (\tau(x)=0^{\scriptscriptstyle{\wedge}} \wedge
 \tau(y)=0^{\scriptscriptstyle{\wedge}} \to
 \tau(x\vee y)=0^{\scriptscriptstyle{\wedge}})]\!]
 \\
 =\bigwedge\limits_{x,y\in
 X}[\![\tau(x^{\scriptscriptstyle{\wedge}})=0^{\scriptscriptstyle{\wedge}}
 \wedge
 \tau(y^{\scriptscriptstyle{\wedge}})=0^{\scriptscriptstyle{\wedge}} \to
 \tau((x\vee y)^{\scriptscriptstyle{\wedge}})=0^{\scriptscriptstyle{\wedge}}]\!]
 =\mathbb{1}.\ \endproof
 \endgathered
 $$

   \subsubsec{3.3.7.}~Theorems 3.3.4 and 3.3.6 were obtained in
\cite{Kut1} and \cite{KutH}, respectively. On using of the above
terminology, the Meyer Theorem (see \cite[3.3.1\,(5)]{DOP} and
\cite{DuM, Mey1}) reads as follows: {\sl Each order bounded
disjointness preserving operator between vector lattices is a
two-point relation.} This fact can be easily deduced from Theorem
3.3.6, since $\ker(bT)$ is a vector sublattice, whenever $T$ is
disjointness preserving.

 \subsection{3.4. Sums of Lattice Homomorphisms}

 In this section we will give a description for an order
 bounded operator $T$ whose modulus may be presented as the sum
 of two lattice homomorphisms in terms of the properties of the
 kernels of the strata of~$T$. Thus, we reveal the connection
 between the $2$-disjoint operators and Grothendieck subspaces.


 \subsubsec{3.4.1.}~Recall that a subspace $H$ of a vector lattice is
 a~{\it $G$-space\/} or {\it Grothendieck subspace\/} provided
 that $H$ enjoys the following property:
 $$
 (\forall\,x,y\in H)\ (x\vee y\vee 0 + x\wedge y\wedge 0 \in H).
 \eqno(1)
 $$

 \subsubsec{3.4.2.}~This condition appears as follows. In 1955
 Grothendieck \cite{Gro} pointed out the subspaces with the above condition
 in the vector lattice $C(Q,\mathbb{R})$ of continuous functions on a compact
 space $Q$ defining them by means of a family of relations $\mathrm{A}$ with
 each relation $\alpha\in\mathrm{A}$ being the form:
 $$
 f(q^1_\alpha)=\lambda_\alpha f(q^2_\alpha)\quad
 (q^1_\alpha,q^2_\alpha\in Q;\,\lambda_\alpha\in\mathbb{R};\,
 \alpha\in\mathrm{A}).
 $$
 These spaces yield examples of $L^1$-predual Banach spaces which are not $AM$-spaces.
 In 1969 Lindenstrauss and Wulpert gave a characterization of such subspaces by means of the
 property 3.4.1 and introduced the term $G$-space (see \cite{LW}). Some related properties of
 Grothendieck spaces are presented also in \cite{Lac} and \cite{Sem}.

 \subsubsec{3.4.3.}~\proclaim{Theorem.}Let $\mathbb{F}$ be a dense subfield in $\mathbb{R}$
 and $X$ a vector lattice over~$\mathbb{F}$. The modulus of an order bounded $\mathbb{F}$-linear
 functional from $X$ to $\mathbb{R}$ is the sum of two lattice homomorphisms if and only if the
 kernel of this functional is a~Grothen\-dieck subspace of $X$.
 \Endproc

 \beginproof~The proof relied on Decomposition Theorem 3.3.4; cp.~\cite{KutG}.~\endproof

 \subsubsec{3.4.4.}~\theorem{}Let $X$ and $Y$ be vector lattices with $Y$ Dedekind complete. The modulus of an order bounded operator $T: X\to Y$ is the sum of some pair of lattice homomorphisms if and only if the kernel of each stratum $bT$ of~$T$ with $b\in\mathbb{B}\!:=\mathbb{P}(Y)$ is a Grothendieck subspace of the ambient vector lattice~$X$.\Endproc

 $\lhd$~The proof runs along the lines of \S~3.3. We apply           the technique of Boolean valued
 ``scalarization'' reducing operator problems to the case of functionals handled in 3.4.3.
 Put $Y=\mathcal{R}{\downarrow}$ and $\tau:=T{\uparrow}$ and the further work with $\tau$ is performed
 within~$\mathbb{V}^{(\mathbb{B})}$. First we observe a useful calculation:
 $$
 [\![\ker(l)\ {\text{is a Grothendieck subspace of}}\
 X^{\scriptscriptstyle{\wedge}} ]\!]
 $$
 $$
 \gathered
 =[\![(\forall\,x,y\in X^{\scriptscriptstyle{\wedge}})
 (\tau(x)=0^{\scriptscriptstyle{\wedge}}
 \wedge
 \tau(y)=0^{\scriptscriptstyle{\wedge}} \to \tau(x\vee y\vee 0+x\wedge
 y\wedge 0)=0^{\scriptscriptstyle{\wedge}})]\!]
 \\
 =\bigwedge\limits_{x,y\in
 X}[\![\tau(x^{\scriptscriptstyle{\wedge}})=0^{\scriptscriptstyle{\wedge}}
 \wedge
 \tau(y^{\scriptscriptstyle{\wedge}})=0^{\scriptscriptstyle{\wedge}}
 \to
 \tau((x\vee y\vee 0+x\wedge y\wedge
 0)^{\scriptscriptstyle{\wedge}})=0^{\scriptscriptstyle{\wedge}}]\!].
 \endgathered
 \eqno(\ast)
 $$

 {\it Sufficiency:} Take $x,y\in X$ and put $b:=[\![Tx=0^{\scriptscriptstyle{\wedge}}]\!]\wedge
 [\![Ty=0^{\scriptscriptstyle{\wedge}}]\!]$. It follows from $(G)$ (see the Gordon Theorem) that  $x,y\in\ker(bT)$.
 By hypothesis $\ker(bT)$ is a Grothendieck subspace and hence $bT(x\vee y\vee 0+x\wedge y\wedge 0)=0$. By using $(G)$
 again we get
 $$
 [\![Tx=0^{\scriptscriptstyle{\wedge}}]\!]\wedge
 [\![Ty=0^{\scriptscriptstyle{\wedge}}]\!]=b\le [\![T(x\vee y\vee
 0+x\wedge y\wedge 0)=0^{\scriptscriptstyle{\wedge}}]\!].
 $$
 Now, it follows from ($\ast$) that $[\![\ker(l)\ {\text{~is a Grothendieck subspace of}}\
 X^{\scriptscriptstyle{\wedge}} ]\!]=\mathbb{1}$. By Transfer Principle we can apply Theorem 3.4.3 to $\tau$ within $\mathbb{V}^{(\mathbb{B})}$, consequently, $|\tau|=\tau_1+\tau_2$ with $\tau_1$ and $\tau_2$ being lattice homomorphisms within $\mathbb{V}^{(\mathbb{B})}$. It can be easily seen that the operators $T_1\!:=\tau_1{\downwardarrow}$
 and $T_2\!:=\tau_2{\downwardarrow}$ from $X$ to $\mathcal{R}{\downarrow}$ are lattice homomorphisms and $|T|=T_1+T_2$.

 {\it Necessity:} Assume that $|T|=T_1+T_2$ for some lattice homomorphisms $T_1,T_2:X\to Y$ and denote $\tau\!:=T{\upwardarrow}$, $\tau_1\!:=T_1{\upwardarrow}$
 and $\tau_2\!:=T_2{\upwardarrow}$. It can be easily checked that inside $\mathbb{V}^{(\mathbb{B})}$ we have $\tau,\tau_1,\tau_2:X^{\scriptscriptstyle{\wedge}}\to\mathcal{R}$
 and $|\tau|=\tau_1+\tau_2$; moreover, $\tau_1$ and $\tau_2$ lattice homomorphisms. By Theorem 3.4.3 and Transfer Principle $[\![\ker(l)\ {\text{ is a Grothendieck subspace of}}\
 X^{\scriptscriptstyle{\wedge}} ]\!]=\mathbb{1}$.
 Making use of ($\ast$), we infer
 $$
 [\![\tau(x^{\scriptscriptstyle{\wedge}})=0^{\scriptscriptstyle{\wedge}}
 \wedge
 \tau(y^{\scriptscriptstyle{\wedge}})=0^{\scriptscriptstyle{\wedge}}]\!]\leq
 [\![\tau((x\vee y\vee 0+x\wedge y\wedge
 0)^{\scriptscriptstyle{\wedge}})=0^{\scriptscriptstyle{\wedge}}]\!].
 $$
 Now, if $b\in\mathbb{B}$ and $bTx=bTy=0$ then
 $[\![l((x\vee y\vee 0+x\wedge y\wedge
 0)^{\scriptscriptstyle{\wedge}})=0^{\scriptscriptstyle{\wedge}}]\!]
 \geq b$, whence by the Gordon Theorem we get
 $bT(x\vee y\vee 0+x\wedge y\wedge 0)=0$.~\endproof

 \subsubsec{3.4.5.}~The main result of the section (Theorem 3.4.4) was obtained in \cite{KutG}. The sums of Riesz homomorphisms were first described in \cite{BHP} in terms of $n$-disjoint operators; see \S~3.3. A survey of some conceptually close results
on $n$-dis\-joint operators is given in~[3, \S\,5.6].

\subsection{3.5. 	Disjointness Preserving Bilinear Operators}
It was observed above that
a~linear operator $T$ from a vector lattice~$X$ to a~Dedekind complete
vector lattice $Y$ is in a~sense determined up to an orthomorphism
from the family of the kernels of the strata $\pi T$ of~$T$
with $\pi$ ranging over all band projections on $Y$. Similar reasoning
was involved in \cite{KK_DPO} to characterize order bounded
disjointness preserving bilinear operators. Unfortunately, Theorem 3.4
in \cite{KK_DPO} is erroneous and this note aims to give the correct
statement and proof of this result.

 In what follows $X$, $Y$, and $Z$ are Archimedean vector lattices, $Z^\mathbb{u}$ is a universal
 completion of $Z$, and $B:X\times Y\rightarrow Z$ is a bilinear operator. We denote  the~Boolean
 algebra of band projections in $X$ by  $\mathbb{P}(X)$. Recall that a linear operator $T:X\rightarrow Y$
 is said to be  \textit{disjointness preserving\/} if $x\perp y$ implies $Tx\perp Ty$ for all $x,y\in X$.
 A bilinear operator $B:X\times Y\rightarrow Z$ is called \textit{disjointness preserving\/}
 (a \textit{lattice bimorphism\/}) if the linear operators $B(x,\cdot):y\mapsto B(x, y)$ $(y\in Y)$
 and $B(\cdot,y):x\mapsto B(x,y)$ $(x\in X)$ are disjointness preserving for all $x\in X$ and $y\in Y$
 (lattice homomorphisms for all $x\in X_+$ and $y\in Y_+$).
 Denote $X_\pi\!:=\bigcap\{\ker(\pi B(\cdot,y)):\ y\in Y\}$ and $Y_\pi\!:=\bigcap\{\ker(\pi B(x,\cdot)):\ x\in X\}$. Clearly, $X_\pi$ and $Y_\pi$ are vector subspaces of $X$ and $Y$, respectively. Now we state the main result of the note.
 \smallskip

 \subsubsec{3.5.1.}~\proclaim{Theorem.}Assume that $X$, $Y$, and $Z$
 are vector lattices with $Z$ having the projection property. For an order
bounded bilinear operator $B:X\times Y\rightarrow Z$ the following  are equivalent:
 \vspace{2pt}

 {\bf(1)}~$B$ is disjointness preserving.
 \vspace{2pt}

 {\bf(2)}~There are a band projection $\varrho\in\mathbb{P}(Z)$ and lattice homomorphisms $S:X\rightarrow Z^\mathbb{u}$ and
 $T:Y\rightarrow Z^\mathbb{u}$ such that
 $B(x,y)=\varrho S(x)T(y)-\varrho^\perp S(x)T(y)$ for all $(x,y)\in X\times Y$.%
 \vspace{2pt}

 {\bf(3)}~For every $\pi\in\mathbb{P}(Z)$ the subspaces $X_\pi$ and $Y_\pi$ are order ideals respectively in $X$ and $Y$, and the kernel of every stratum $\pi B$ of $B$ with $\pi\in\mathbb{P}(Z)$ is representable as
 $$
 \ker(\pi B)=\bigcup\big\{X_{\sigma}\times Y_{\tau}:\ \sigma,\tau\in\mathbb{P}(Z);\,\sigma\vee\tau=\pi\big\}.
 $$
 \Endproc

 The proof  proceeds  along the general lines of \cite{KutH}--\cite{KK_DPO}: Using the canonical embedding and ascent to the Boolean valued universe $\mathbb{V}^{(\mathbb{B})}$, we reduce the matter to characterizing a disjointness preserving bilinear functional on the product of two vector lattices over a dense subfield of the reals $\mathbb{R}$. The resulting scalar problem is solved by the following simple fact.
 \smallskip

 \subsubsec{3.5.2}\proclaim{}Let $X$ and $Y$ be vector lattices. For an order bounded bilinear functional $\beta:X\times Y\rightarrow\mathbb{R}$ the following  are equivalent:
 \vspace{2pt}

 {\bf (1)}~$\beta$ is disjointness preserving.
 \vspace{2pt}

 {\bf(2)}~$\ker(\beta)=(X_0\times Y)\cup(X\times Y_0)$ for some order ideals $X_0\subset X$ and $Y_0\subset Y$.
 \vspace{2pt}

 {\bf(3)}~There exist lattice homomorphisms $g:X\rightarrow\mathbb{R}$ and $h:Y\rightarrow\mathbb{R}$ such that either $\beta(x,y)=g(x)h(y)$ or $\beta(x,y)=-g(x)h(y)$ for all $x\in X$ and $y\in Y$.
 \Endproc

 \beginproof~Assume that $\ker(\beta)=(X_0\times Y)\cup(X\times Y_0)$ and take $y\in Y$. If $y\in Y_0$ then $\beta(\cdot,y)\equiv0$, otherwise $\ker(\beta(\cdot,y))=X_0$ and $\beta(\cdot,y)$ is disjointness preserving, since an order bounded linear functional is disjointness preserving if and only if its null-space is an order ideal. Similarly,  $\beta(x,\cdot)$ is disjointness preserving for all $x\in X$ and thus $(2)\Longrightarrow(1)$. The implication $(1)\Longrightarrow(3)$ was established in \cite[Theorem 3.2]{KT1} and $(3)\Longrightarrow(1)$ is trivial with $X_0=\ker(g)$ and $Y_0=\ker(h)$.~\endproof
 \smallskip

 Let $\mathbb{B}$ be a complete Boolean algebra and let $\mathbb{V}^{(\mathbb{B})}$ be the corresponding Boolean valued model with the Boolean truth value $[\![\varphi]\!]$ of a set-theoretic formula $\varphi$. There exists $\mathcal{R}\in\mathbb{V}^{(\mathbb{B})}$ playing the role of the~field of reals within $\mathbb{V}^{(\mathbb{B})}$. The descent functor sends each internal algebraic structure $\mathfrak{A}$ into its descent $\mathfrak{A}{\downarrow}$ which is an algebraic structure in the conventional sense. Gordon's theorem (see \cite[8.1.2]{DOP} and \cite[Theorem 2.4.2]{KK_T}) tells us that the algebraic structure $\mathcal{R}{\downarrow}$ $($with the descended operations and order$)$ is a universally complete vector lattice. Moreover, there is a Boolean isomorphism $\chi$ of $\mathbb{B}$ onto $\mathbb{P}(\mathcal{R}{\downarrow})$ such that
 $b\le[\![\,x=y\,]\!]$ if and only if $\chi(b)x=\chi(b)y$. We identify $\mathbb{B}$ with $\mathbb{P}(\mathcal{R}{\downarrow})$ and take $\chi$ to be $I_\mathbb{B}$.

 Let $[X\times Y,\mathcal{R}{\downarrow}]\in\mathbb{V}$ and let
 $[\![X^{\scriptscriptstyle\wedge}\times Y^{\scriptscriptstyle\wedge}, \mathcal{R}]\!]\in\mathbb{V}^{(\mathbb{B})}$
 stand respectively for the sets  of all maps from $X\times Y$ to $\mathcal{R}{\downarrow}$ and from $X^{\scriptscriptstyle\wedge}\times X^{\scriptscriptstyle\wedge}$ to $\mathcal{R}$ (within $\mathbb{V}^{(\mathbb{B})}$). The correspondences $f\mapsto f{\upwardarrow}$, the modified ascent, is a bijection between $[X\times Y,\mathcal{R}{\downarrow}]$ and
 $[\![X^{\scriptscriptstyle\wedge}\times Y^{\scriptscriptstyle\wedge}, \mathcal{R}\,]\!]$.
 Given $f\in[X,\mathcal{R}{\downarrow}]$, the internal map $f{\upwardarrow}\in[\![X^{\scriptscriptstyle\wedge},\mathcal{R}]\!]$ is uniquely determined by the relation $[\![f{\upwardarrow}(x^{\scriptscriptstyle\wedge})=f(x)]\!] =\mathbb{1}$ $(x\in X)$. Observe also that $\pi\leq [\![f{\upwardarrow}(x^{\scriptscriptstyle\wedge})=\pi f(x)]\!]$ $(x\in X,\,\pi\in\mathbb{P}(\mathcal{R}{\downarrow}))$. This fact specifies for bilinear operators as follows:

  \smallskip

 \subsubsec{3.5.3.}\proclaim{}Let $B:X\times Y\rightarrow Y$ be a bilinear operator and $\beta\!:=B{\upwardarrow}$ its modified ascent. Then $\beta:X^{\scriptscriptstyle\wedge}\times Y^{\scriptscriptstyle\wedge}\rightarrow\mathcal{R}$ is a $\mathbb{R}^{\scriptscriptstyle\wedge}$-bilinear functional within $\mathbb{V}^{(\mathbb{B})}$. Moreover, $B$ is order bounded and disjointness preserving if and only $[\![\,\beta$ is order bounded and disjointness preserving $]\!]=\mathbb{1}$.
 \Endproc

 \beginproof~The proof goes along similar lines to the proof of Theorem 3.3.3 in \cite{KK_T}.~\endproof
 \smallskip

  \subsubsec{3.5.4.}\proclaim{}Let $B$ and $\beta$ be as in 3.4.3. Then
 $[\![\ker(B)^{\scriptscriptstyle\wedge}=\ker(\beta)]\!]=\mathbb{1}$.
  \Endproc

 \beginproof~Using the above-mentioned  property of the modified ascent and interpreting the formal definition $z\in\ker(\beta)\leftrightarrow(\exists x\in X^{\scriptscriptstyle\wedge})(\exists y\in Y^{\scriptscriptstyle\wedge})(z=(x,y)\wedge\beta(x,y)=0)$, the proof  reduces to the straightforward calculation:
 \begin{align*}
 [\![z\in\ker(\beta)]\!]&=\bigvee_{x\in X,\,y\in Y}
 [\![z=(x^{\scriptscriptstyle\wedge},y^{\scriptscriptstyle\wedge})\wedge
 \beta(x^{\scriptscriptstyle\wedge},y^{\scriptscriptstyle\wedge})=0]\!]
 \\
 &=\bigvee_{(x,y)\in X\times Y}
 [\![z=(x,y)^{\scriptscriptstyle\wedge}\wedge (x,y)^{\scriptscriptstyle\wedge}\in\ker(B)^{\scriptscriptstyle\wedge}]\!]
 \end{align*}
 \begin{align*}
 \leq[\![z\in\ker(B)^{\scriptscriptstyle\wedge}]\!]
 &=\bigvee_{(x,y)\in X\times Y} [\![z=(x,y)^{\scriptscriptstyle\wedge}\wedge (x.y)\in\ker(B)]\!]
 \\
 &=\bigvee_{x\in X,\,y\in Y}
 [\![(z=(x^{\scriptscriptstyle\wedge},y^{\scriptscriptstyle\wedge})\wedge
 \beta(x^{\scriptscriptstyle\wedge},y^{\scriptscriptstyle\wedge})=0]\!]
 \\
 \leq[\![z\in\ker(\beta)]\!].\ &\ \endproof
 \end{align*}

 \subsubsec{3.4.5.}~\proclaim{}Define $\mathcal{X}$ and $\mathcal{Y}$ within $\mathbb{V}^{(\mathbb{B})}$ as follows: $\mathcal{X}\!:=\bigcap\{\ker(\beta(\cdot,Y)):\,y\in Y^{\scriptscriptstyle\wedge}\}$ and $\mathcal{Y}\!:=\bigcap\{\ker(\beta(x,\cdot)):\,x\in X^{\scriptscriptstyle\wedge}\}$. Given arbitrary $\pi\in\mathbb{P}(Z)$, $x\in X$, and $y\in Y$, we have
 $$
 \pi\leq[\![x^{\scriptscriptstyle\wedge}\in\mathcal{X}]\!] \Longleftrightarrow x\in X_\pi,\quad \pi\leq[\![y^{\scriptscriptstyle\wedge}\in\mathcal{Y}]\!] \Longleftrightarrow y\in Y_\pi.%
 $$
 \Endproc

 \beginproof~Given $\pi\in\mathbb{P}(Z)$ and $x\in X$, we need only to calculate Boolean truth values taking it into account that $[\![B(x,y)=\beta(x^{\scriptscriptstyle\wedge},v^{\scriptscriptstyle\wedge})]\!]=\mathbb{1}$
 for all $x\in X$ and $y\in Y$:
 $$
 [\![x^{\scriptscriptstyle\wedge}\in\mathcal{X}]\!]
 =[\![(\forall v\in Y^{\scriptscriptstyle\wedge})
 \beta(x^{\scriptscriptstyle\wedge},v)=0]\!]
 \\
 =\bigwedge_{v\in Y} [\![\beta(x^{\scriptscriptstyle\wedge},v^{\scriptscriptstyle\wedge})=0]\!]
 =\bigwedge_{v\in Y} [\![B(x,v)=0]\!]
 $$
 It follows that $\pi\leq[\![x^{\scriptscriptstyle\wedge}\in\mathcal{X}]\!]$ if and only if $\pi\leq[\![B(x,v)=0]\!]$ for all $v\in Y$. By Gordon's theorem the latter means that $\pi B(x,v)=0$ for all $v\in Y$, that is $x\in X_\pi$.~\endproof
 \smallskip

 \subsubsec{3.4.6.}\proclaim{}Let $B$ and $\beta$ be as in 3.4.3. For all $\pi\in\mathbb{P}(Z)$, $x\in X$, and $y\in Y$, we have $\pi\leq[\![(x^{\scriptscriptstyle\wedge},y^{\scriptscriptstyle\wedge})\in
 (\mathcal{X}\times Y)\cup(X\times\mathcal{Y})]\!]$
 if and only if there exist $\sigma,\tau\in\mathbb{P}(Z)$ such that $\sigma\vee\tau=\pi$, $x\in X_\sigma$, and $y\in Y_\tau$.
 \Endproc

 \beginproof~Put $\rho\!:=[\![(x^{\scriptscriptstyle\wedge},y^{\scriptscriptstyle\wedge})\in
 (\mathcal{X}\times Y)\cup(X\times\mathcal{Y})]\!]$ and observe that
 $$
 \rho=[\![(x^{\scriptscriptstyle\wedge}\in\mathcal{X})\vee y^{\scriptscriptstyle\wedge}\in\mathcal{Y}]\!]
 =[\![x^{\scriptscriptstyle\wedge}\in\mathcal{X}]\!]\vee [\![y^{\scriptscriptstyle\wedge}\in\mathcal{Y}]\!].
 $$
 Clearly, $\pi\leq\rho$ if and only if $\sigma\vee\tau=\pi$ for some $\sigma\leq[\![x^{\scriptscriptstyle\wedge}\in\mathcal{X}]\!]$ and $\tau\leq[\![y^{\scriptscriptstyle\wedge}\in\mathcal{Y}]\!]$, so that the required property follows from
 3.4.5.~\endproof
 \smallskip

 {\sc Proof of Theorem 3.5.1.}~The implication $(1)\Longrightarrow(2)$ was proved in \cite[Corollary~3.3]{KT1}, while $(2)\Longrightarrow(3)$ is straightforward. Indeed, observe first that if (2) is fulfilled then $|B(x,y)|=|B|(|x|,|y|)=|S|(|x|)|T|(|y|)$, so that we can assume $S$ and $T$ to be lattice homomorphisms, as in this event $\ker(B)=\ker(|B|)$. Take $\pi\in\mathbb{P}(Z)$ and put $\sigma\!:=\pi-\pi[Sx]$ and $\tau\!:=\pi-\pi[Ty]$, where $[y]$ is a band projection onto $\{y\}^{\perp\perp}$. Observe next that $\pi B(x,y)=0$ if and only if $\pi[Sx]$ and $\pi[Ty]$ are disjoint or, which is the same, if $\sigma\vee\tau=\pi$. Moreover, the map $\rho_y:x\mapsto\sigma S(x)T(y)$ is disjointness preserving for all $y\in Y$ and so $X_\sigma=\bigcap_{y\in Y}\ker(\rho_y)$ is an order ideal in $X$. Similarly, $Y_\tau$ is an order ideal in $Y$. Thus, $(x,y)\in\ker(\pi B)$ if and only if
 $x\in X_\sigma$ and $y\in Y_\tau$ for some $\sigma,\tau\in\mathbb{P}(Z)$ with $\sigma\vee\tau=\pi$.

 Prove the remaining implication $(3)\Longrightarrow(1)$. Suppose that
(3) holds for all $\pi\in\mathbb{P}(Y)$. Take $x,u\in X$ and put $\pi\!:=[\![x^{\scriptscriptstyle\wedge}\in\mathcal{X}]\!]$ and $\rho\!:=[\![|u|^{\scriptscriptstyle\wedge}\leq |x|^{\scriptscriptstyle\wedge}]\!]$. By 3.4.5 we have $x\in X_\pi$. Note also that either $\rho=\mathbb{0}$ or $\rho=\mathbb{1}$. If $\rho=\mathbb{1}$ then $|u|\leq|x|$ and by hypotheses $u\in X_\pi$. Again by 3.4.6 we get $\rho\leq[\![u^{\scriptscriptstyle\wedge}\in\mathcal{X}]\!]$. This estimate is obvious whenever $\rho=\mathbb{0}$, so that $[\![x^{\scriptscriptstyle\wedge}\in\mathcal{X}]\!]\wedge
 [\![|u|^{\scriptscriptstyle\wedge}\leq |x|^{\scriptscriptstyle\wedge}]\!]\Rightarrow
 [\![u^{\scriptscriptstyle\wedge}\in\mathcal{X}]\!]=\mathbb{1}$ for all $x,u\in X$. Simple calculation shows
 now that $\mathcal{X}$ is an order ideal in $X^{\scriptscriptstyle\wedge}$:
 \begin{multline*}
 [\![(\forall x,u\in X^{\scriptscriptstyle\wedge})(|u|\leq|x|\,\wedge\, x\in\mathcal{X}\rightarrow u\in\mathcal{X})]\!]
 \\
 =\bigwedge_{u,x\in X}\big([\![x\in\mathcal{X}]\!] \wedge[\![|u|\leq|x|]\!]\Rightarrow[\![u\in\mathcal{X}]\!]\big)=\mathbb{1}.
 \end{multline*}
 Similarly, $\mathcal{Y}$ is an order ideal in $Y^{\scriptscriptstyle\wedge}$.

 It follows from  (3) and 3.4.6 that $(x,y)\in \ker(\pi B)$ if and only if $\pi\leq[\![(x^{\scriptscriptstyle\wedge},y^{\scriptscriptstyle\wedge})\in
 (\mathcal{X}\times Y)\cup(X\times\mathcal{Y})]\!]$. Taking into account 3.4.3 and the observation  before it we conclude that $\pi\leq[\![(x^{\scriptscriptstyle\wedge}, y^{\scriptscriptstyle\wedge})\in\ker(\beta)]\!])$ if and only if $\pi\leq[\![(x^{\scriptscriptstyle\wedge}, y^{\scriptscriptstyle\wedge})\in(\mathcal{X}\times Y)\cup(X\times\mathcal{Y})]\!]$ and so $[\![\ker(\beta)=
 (\mathcal{X}\times Y)\cup(X\times\mathcal{Y})]\!] =\mathbb{1}$. It remains to apply the equivalence (1)$\Longleftrightarrow$(3) of 3.4.2 within $\mathbb{V}^{(\mathbb{B})}$ . It follows that $B$ is disjointness preserving by 3.4.3.~\endproof

 \smallskip

 \subsubsec{3.4.7.}\proclaim{Corollary.}Assume that $Y$ has the projection property. An order bounded linear operator $T:X\rightarrow Y$ is disjointness preserving if and only if $\ker(bT)$ is an order
 ideal in $X$ for every projection $b\in\mathbb{P}(Y)$.
 \Endproc

 \beginproof~Apply the above theorem to the bilinear operator $B:X\times\mathbb{R}\rightarrow Y$ defined as $B(x,\lambda)=\lambda T(x)$ for all $x\in X$ and $\lambda\in\mathbb{R}$.~\endproof

 \section{Chapter 4. Order Continuous Operators}

 \subsection{4.1.~Maharam Operators}

 Now we examine some class of order continuous positive operators that behave in many instances like functionals. In fact such operators are representable as Boolean valued order continuous functionals.

 \subsubsec{4.1.1.}~Throughout this section $X$ and $Y$ are vector lattices with $Y$ Dedekind complete. A linear operator
 $T: X\to Y$ is said to have the {\it Maharam
 property\/} or is said to be \textit{order interval preserving} whenever $T[0,x]=[0, Tx]$ for every $0\le x\in X$; i.~e., if for every $0\le x\in X$ and $0\le y\le Tx$ there is
 some~$0\le u\in X$ such that $Tu=y$ and $0 \le u\le x$. A~{\it Maharam operator} is an order continuous positive operator whose modulus enjoys the Maharam property.

 Say that a linear operator $S:X\to Y$ is \textit{absolutely
 continuous\/} with respect to $T$ and write $S\preccurlyeq T$
 if $|S|x\in\{|T|x\}^{\perp\perp}$ for all $x\in X_+$. It can
 be easily seen that if $S\in\{T\}^{\perp\perp}$ then
 $S\preccurlyeq T$, but the converse may be false.

 \subsubsec{4.1.2.}~The null ideal $\mathcal{N}_T$ of an order
 bounded operator $T:X\rightarrow Y$ is defined by
 $\mathcal{N}_T\!:=\{x\in X:\,|T|(|x|)=0$. Observe that
 $\mathcal{N}_T$ is indeed an ideal in $X$. The disjoint
 complement of $\mathcal{N}_T$ is referred to as the \textit{carrier} of $T$ and
 is denoted by $\mathcal{C}_T$, so that $\mathcal{C}_T\!:=\mathcal{N}_T^{\perp}$.
 An operator $T$ is called \textit{strictly positive\/} whenever $0<x\in X$ implies $0<|T|(x)$.
 Clearly, $|T|$ is strictly positive on $\mathcal{C}_T$. Sometimes we find it convenient to~denote
 $X_T\!:=\mathcal{C}_T$ and $Y_T\!:=(\im T)^{\perp\perp}$.

 \subsubsec{4.1.3.}~As an examples of Maharam operators, we consider conditional expectation
 and Bochner integration. Take a~probability space $(Q,\Sigma,\mu)$ and let $\Sigma_0$ and
 $\mu_0$ be a~$\sigma$-subalgebra of $\Sigma$ and the restriction of $\mu$ to $\Sigma_0$.
 The conditional expectation operator $\mathcal{E}(\cdot,\Sigma_0)$ is a~Maharam operator from
 $L^1(Q,\Sigma,\mu)$ onto $L^1(Q,\Sigma_0,\mu_0)$.The restriction of $\mathcal{E}(\cdot,\Sigma_0)$
 to $L^p(Q,\Sigma,\mu)$ is also a~Maharam operator from $L^p(Q,\Sigma,\mu)$ to $L^p(Q,\Sigma_0,\mu_0)$.
 These facts are immediate in view of  the simple properties of conditional expectation.

 Let $(Q,\Sigma,\mu)$ be a~probability space, and let $Y$ be a~Banach lattice. Consider
 the space $X\!:=L^1\,(Q,\Sigma,\mu,F)$ of Bochner integrable $Y$-valued functions, and
 let $T:E\rightarrow F$ denote the Bochner integral $Tf\!:=\int_Qf\,d\mu$. If the Banach
 lattice $Y$ is has order continuous norm (in this case $Y$ is order complete) then $X$
 is a Dedekind complete vector lattice under the natural order $f\geq0\Leftrightarrow f(t)\geq0$ for almost all
 $t\in Q$) and $T$ is a~Maharam operator.  See more examples in \cite{DOP, DOR}.

 \subsubsec{4.1.4.}~A positive operator $T:X\to Y$ is said to have the
 \textit{Levi property\/} if~$\sup x_\alpha$
 exists in $X$ for every increasing net $(x_\alpha)\subset X_+$,
 provided that the net $(Tx_\alpha)$ is order bounded in~$Y$.
 For an order bounded order continuous operator $T$ from $X$
 to $Y$ denote by $\mathcal{D}_m(T)$ the largest ideal of the
 universal completion $X^{\rm u}$ onto which we may extend
 the operator $T$ by order continuity. For a positive order
 continuous operator $T$ we have $X=\mathcal{D}_m(T)$ if and
 only if $T$ has the Levi property.

 The following theorem describes an important property of Maharam operators, enabling us
 to embed them into an appropriate Boolean-valued universe as order continuous functionals.

 \subsubsec{4.1.5.}~\theorem{}Let $X$ and $Y$ be some vector
 lattices with $Y$ having the projection property and let $T$
 be a~Maharam operator from $X$ to $Y$. Then there exist an
 order closed subalgebra\/ $\mathcal{B}$ of\/ $\mathbb{B}(X_T)$
 consisting of projection bands and a~Boolean isomorphism $h$
 from $\mathbb{B}(Y_T)$ onto $\mathcal{B}$ such that
 $T(h(L))\subset L$ for all $L\in\mathbb{B}(Y_T)$.
 \Endproc

 The Boolean algebra of projections $\mathcal{B}$ in Theorem 4.1.5 as well as the corresponding Boolean algebra of bands admits a simple description. For $L\in\mathbb{B}(Y_T)$ denote by $h(L)$ the band in $\mathbb{B}(X_T)$ corresponding to the band projection $h([K])$.

 \subsubsec{4.1.6.}\proclaim{}~For a band $K\in\mathbb{B}(X_T)$ the following are equivalent:

 \subsubsec{(1)}~$Tu=Tv$ and $u\in K$ imply $v\in K$ for all $u,v\in X_+$.

 \subsubsec{(2)}~$T(K'_+)\subset T(K_+)^{\perp\perp}$ implies $K'\subset K$
 for all $K'\in\mathbb{B}(X_T)$.

 \subsubsec{(3)}~$K=h(L)$ for some $L\in\mathbb{B}(Y)$.
 \Endproc

 A band $K\in\mathbb{B}(X_T)$ (as well as the corresponding
 band projection $[K]\in\mathbb{P}(X_T)$) is said to be
 $T$-\textit{saturated\/} if one of (and then all) the conditions 4.1.6\,(1--3) is fulfilled.

 The following can be deduce from 4.1.6 by the Freudenthal Spectral Theorem.%

 \subsubsec{4.1.7.}~\proclaim{}If $X$ and $Y$ are Dedekind complete vector lattices and $T$ is a Maharam operator from $X$ to $Y$, then there exists an $f$-module structure on $X$ over an $f$-algebra $\mathcal{Z}(Y)$ such that an order bounded operator $S:X\to Y$ is absolutely continuous with respect to  $T$ if and only if $S$ is $\mathcal{Z}(Y)$-linear.%
 \Endproc

 We now state the main result of the section.

 \subsubsec{4.1.8.} \theorem{}Let $X$ be a vector lattice,
 $Y\!:=\mathcal{R}{\downarrow}$, and let $T:X\to Y$ be a~strictly positive Maharam operator with $Y=Y_T$. Then there are $\mathcal{X}, \tau\in\mathbb{V}^{(\mathbb{B})}$
 satisfying the following:

 \subsubsec{(1)}~$\mathbb{V}^{(\mathbb{B})}\models\text{~``}\mathcal{X}$
 is a~Dedekind complete vector lattice and
 $\tau:\mathcal{X}\to\mathcal{R}$ is an strictly positive order
 continuous functional with the Levi property''.

 \subsubsec{(2)}~$\mathcal{X}{\downarrow}$ is a~Dedekind complete vector lattice and a unitary $f$-module
 over the $f$-algebra $\mathcal{R}{\downarrow}$.

 \subsubsec{(3)}~$\tau{\downarrow}:
 \mathcal{X}{\downarrow}\to\mathcal{R}{\downarrow}$ is a~strictly positive Maharam operator with the Levi
 property and an $\mathcal{R}{\downarrow}$-module homomorphism.

 \subsubsec{(4)}~There exists a~lattice isomorphism $\varphi$
 from $X$ into $\mathcal{X}{\downarrow}$ such that
 $\varphi(X)\subset X^\delta\subset\mathcal{X}{\downarrow}
 \subset X^{\rm u}$ is a universal completion of $X$ and
 $T={\tau{\downarrow}\circ\varphi}$.
 \Endproc

 \subsubsec{4.1.9.}~The Maharam operators stems from from the theory of Maharam's ``full-valued'' integrals,
 see \cite{Mah3, Mah4, Mah5}. Theorem 4.1.8 was established in \cite{K12, Kus1}. More results, applications, and references
 on Maharam operators are in~\cite{DOP, DOR}.
 See \cite{KK3} for some extension of this theory to sublinear and convex operators.

 \subsection{4.2. Representation of Order Continuous Operators}

 Theorem 4.1.8 enables us to state that
 every fact on order continuous functionals ought to have a~parallel variant for a~Maharam operators which
 may be proved by the Boolean valued machinery. The aim of this section is to prove an operator version of the following result.

 \subsubsec{4.2.1.}~\proclaim{Theorem.}Let $X$ be a vector lattice and
 $X_n^\sim$ separates the points of $X$. Then there exist
 order dense ideals $L$ and $X'$ in $X^{\mathrm{u}}$ and a
 linear functional $\tau:L\to\mathbb{R}$ such that
 \vspace{2pt}

 \subsubsec{(1)}~$X'=\{x'\in X':\
 xx'\in L\text{~for all~}x\in X\}$.
 \vspace{2pt}

 \subsubsec{(2)}~$\tau$ is strictly positive, $o$-continuous,
 and has the Levi property.%
 \vspace{2pt}

 \subsubsec{(3)}~For every $\sigma\in X^\sim_n$ there exists a
 unique $x'\in X'$ such that
 $$
 \sigma(x)=\tau(x\cdot x')\quad(x\in X).
 $$

 \subsubsec{(4)}~The map $\sigma\mapsto x'$ is a
 lattice isomorphism of $X_n^\sim$ onto $X'$.
 \Endproc

 \subsubsec{4.2.2.}~To translate Theorem 4.10.1 into a result on operators we need some preparations. Let $X$ and $Y$ be
 $f$-modules over an $f$-algebra $A$. A linear operator
 $T:X\to Y$ is called $A$-\textit{linear\/} if $T(ax)=aTx$ for
 all $x\in X$ and $a\in A$. Denote by $L^A(X,Y)$ the set of all
 order bounded $A$-linear operators from $X$ to $Y$ and put
 $L^A_n(X,Y)\!:=L^A(X,Y)\cap L^\sim_n(X,Y)$.

 Say that a set $\mathcal{T}\subset L^\sim(X,Y)$
 \textit{separates the points\/} of $X$ whenever, given nonzero
 $x\in X$, there exists $T\in\mathcal{T}$ such that $Tx\ne0$. In the
 case of a~Dedekind complete~$Y$ and the sublattice $\mathcal{T}\subset L^\sim(X,Y)$
 this is equivalent to saying that for every nonzero $x\in X_+$
 there is a positive operator $T\in\mathcal{T}$ with $Tx\ne0$.

 \subsubsec{4.2.3.}~Given a real vector lattice
 $\mathcal{X}$ within $\mathbb{V}^{(\mathbb{B})}$, denote by
 $\mathcal{X}^\sim$ and $\mathcal{X}^\sim_n$ the internal
 vector lattices of order bounded and order continuous
 functionals on~$\mathcal{X}$, respectively. More precisely, $[\![\sigma\in\mathcal{X}^\sim]\!]=\mathbb{1}$ and
 $[\![\sigma\in\mathcal{X}_n^\sim]\!]=\mathbb{1}$ mean that
 $[\![\sigma:\mathcal{X}\to\mathcal{R}$ is an order bounded functional $]\!]=\mathbb{1}$ and $[\![\sigma:\mathcal{X}\to\mathcal{R}$ is an
 order continuous functional $]\!]=\mathbb{1}$, respectively.
 Put $X\!:=\mathcal{X}{\downarrow}$ and
 $A\!:=\mathcal{R}{\downarrow}$.

 \subsubsec{4.2.4.}~\proclaim{Theorem.}~The
 mapping assigning to each
 $\sigma\in\mathcal{X}^\sim{\downarrow}$
 its descent $S\!:=\sigma{\downarrow}$ is a lattice isomorphism
 of $\mathcal{X}^\sim{\downarrow}$ and
 $\mathcal{X}^\sim_n{\downarrow}$ onto
 $L^A(X,\mathcal{R}{\downarrow})$ and
 $L^A_n(X,\mathcal{R}{\downarrow})$, respectively. Moreover,
 $$
 [\![\mathcal{X}^\sim\ (\text{resp.~}\mathcal{X}^\sim_n)
 \text{~separates the points of~}\mathcal{X}]\!]=\mathbb{1}
 $$
 if and only if $L^A(X,\mathcal{R}{\downarrow})$ $($resp.
 $L^A_n(X,\mathcal{R}{\downarrow}))$ separates the points of
 $X$.
 \Endproc

 \subsubsec{4.2.5}~By the Gordon Theorem
 we may assume also that $Y^\mathrm{u}=\mathcal{R}{\downarrow}$. Of~course, in this event we can
 identify $A^\mathrm{u}$ with $Y^\mathrm{u}$. In view of Theorem 1.3.7 there exists a real
 Dedekind complete vector lattice $\mathcal{X}$ within
 $\mathbb{V}^{(\mathbb{B})}$ with $\mathbb{B}=\mathbb{P}(Y)$
 such that $\mathcal{X}{\downarrow}$ is an $f$-module over
 $A^{\mathrm{u}}$, and there is an $f$-module isomorphism $h$
 from $X$ to $\mathcal{X}{\downarrow}$ satisfying
 $\mathcal{X}{\downarrow}=\mix(h(X))$. In virtue of 4.2.4
 $\mathcal{X}_n^\sim$ separates the points of $\mathcal{X}$.
 The Transfer Principle tells us that Theorem 4.2.1 is true
 within $\mathbb{V}^{(\mathbb{B})}$, so that there exist an
 order dense ideal $\mathcal{L}$ in $\mathcal{X}^{\mathrm{u}}$
 and a strictly positive linear functional
 $\tau:\mathcal{L}\to\mathcal{R}$ with the Levi property such
 that the order ideal
 $\mathcal{X}'=\{x'\in\mathcal{X}^{\mathrm{u}}:\
 x'\mathcal{X}\subset\mathcal{L}\}$ is lattice isomorphic to
 $\mathcal{X}^\sim_n$; moreover, the isomorphism is implemented
 by assigning the functional $\sigma_{x'}\in\mathcal{X}^\sim_n$
 to $x'\in\mathcal{X}'$ by $\sigma_{x'}(x)=\tau(xx')$
 $(x\in\mathcal{X})$.

 \subsubsec{4.2.6}~Put $\hat{X}\!:=\mathcal{X}{\downarrow}$,
 $\hat{L}\!:=\mathcal{L}{\downarrow}$,
 $\hat{T}\!:=\tau{\downarrow}$, and
 $\hat{X}'\!:=\mathcal{X}'{\downarrow}$. By~Theorem~2.11.9 we
 can identify the universally complete vector lattices
 $X^{\mathrm{u}}$, $\hat{X}^{\mathrm{u}}$, and
 $\mathcal{X}^{\mathrm{u}}{\downarrow}$ as well as $X$ with a
 laterally dense sublattice in $\hat{X}$. Then $\hat{L}$ is an
 order dense ideal in $\hat{X}^{\mathrm{u}}$ and an $f$-module
 over $A^{\mathrm{u}}$, while $\hat{T}:\hat{L}\to
 Y^{\mathrm{u}}$ is a strictly positive Maharam operator with
 the Levi property. Since the multiplication in
 $X^{\mathrm{u}}$ is the descent of the internal multiplication
 in $\mathcal{X}^{\mathrm{u}}$, we have the representation
 $\hat{X}'=\{x'\in X^{\mathrm{u}}:\,x'\hat{X}\subset\hat{L}\}$.
 Moreover, $\hat{X}'$ is $f$-module isomorphic to
 $L^A_n(\hat{X},Y^{\mathrm{u}})$ by assigning to $x'\in\hat{X}$
 the operator $\hat{S}_{x'}\in L^A_n(\hat{X},Y^{\mathrm{u}})$
 defined as $\hat{S}_{x'}(x)=\hat{T}(xx')$ $(x\in\hat{X})$.
 Now, defining
 $$
 \gathered
 L\!:=\{x\in\hat{L}:\,\hat{T}x\in Y\},\quad T\!:=\hat{T}|_L,
 \\
 X'\!:=\{x'\in\hat{X}':\,x'X\subset L\},
 \endgathered
 $$
 yields that if $x'\in X'$ then $S_{x'}\!:=\hat{S}_{x'}|_X$
 is contained in $L^A_n(X,Y)$. Conversely, an arbitrary
 $S\in L^A_n(X,Y)$ has a representation $Sx=\hat{T}(xx')$
 $(x\in X)$ with some $x'\in\hat{X}'$, so that
 $\hat{T}(xx')\in Y$ for all $x\in X$ and hence $x'\in X'$,
 $xx'\in L$ for all $x\in X$, and $Sx=T(xx')$ $(x\in X)$ by the
 above definitions.

 \subsubsec{4.2.7.}~\proclaim{Theorem.}Let $X$ be an $f$-module over $A\!:=\mathcal{Z}(Y)$ with $Y$ being a~Dedekind complete
 vector lattice and let $L^A_n(X,Y)$ separates the points of
 $X$. Then there exist an order dense ideal $L$ in
 $X^{\mathrm{u}}$ and a strictly positive Maharam operator
 $T:L\to Y$ such that the order ideal $X'=\{x'\in X':\ (\forall\,
 x\in X)\,xx'\in L\}\subset X^\mathrm{u}$ is lattice isomorphic
 to $L^A_n(X,Y)$. The isomorphism is implemented by assigning
 the operator $S_{x'}\in L^A_n(X,Y)$ to an~element $x'\in X'$
 by the formula
 $$
 S_{x'}\,(x)=\Phi\,(xx')\quad(x\in X).
 $$
 If there exists a strictly positive $T_0\in L^A_n(X,Y)$ then one
 can choose $L$ and $T$ such that $X\subset$ L and $T|_X=T_0$.
 \Endproc

 Below, in 4.2.8--4.2.10, $X$ and $Y$ are Dedekind complete vector lattices.

 \subsubsec{4.2.7.} \proclaim{Hahn Decomposition Theorem.}~Let
 $S :X\to Y$ be a~Maharam operator. Then there is a~band
 projection $\pi\in\mathbb{P}(X)$ such that $S^+ =S\circ\pi$
 and $S^-=-S\circ\pi^\perp $. In particular,
 $|S|=S\circ(\pi-\pi^\perp)$.
 \Endproc%

 \subsubsec{4.2.8.}~\proclaim{Nakano Theorem.}~Let
 $T_1,T_2:X\to Y$
 be order bounded operators such that $T\!:=|T_1|+|T_2|$ is
 a~Maharam operator. Then $T_1$ and $T_2$ are disjoint if and
 only if so are their carriers; symbolically,
 $T_1\perp T_2\Longleftrightarrow \mathcal{C}_{T_1}\perp\mathcal{C}_{T_2}$.
 \Endproc

 \subsubsec{4.2.9.}~\proclaim{Radon--Nikod\'ym Theorem.}~Assume that $T:X\to Y$ be a positive Maharam operator. A positive operator $S:X\to Y$ belongs to $\{T\}^{\perp\perp}$ if and only if there exists an orthomorphism $0\le\rho\in\Orth^\infty(X)$ with $Sx=T(\rho x)$ for all $x\in\mathcal{D}(\rho )$.
 \Endproc

 \subsubsec{4.2.11.} Theorem 4.2.1 is proved in~\cite[Theorem 2.1]{VL}. It can be also extracted from \cite[Theorem IX.3.1]{Vul} or \cite[Theorem 3.4.8]{DOP}. Theorem 4.2.7 was proved in \cite{K12}; also see~\cite{DOP}. Theorems 4.2.8--4.2.10, first obtained in \cite{LS}, may be deduced easily from 4.2.7, or can be proved by the general scheme of ``Boolean valued scalarization.''

 \subsection{4.3. Conditional Expectation Type Operators}

 The conditional expectation operators have many remarkable
 properties related to the order structure of the underlying
 function space. Boolean valued analysis enables us to
 demonstrate that some much more general class of operators shares these properties.

 \subsubsec{4.3.1.}~Let $Z$ be a universally complete vector lattice with
 unit $\mathbb{1}$. Recall that $Z$ is an $f$-algebra with multiplicative
 unit $\mathbb{1}$. Assume that $\Phi:L^1(\Phi)\to Y$ is a Maharam operator with
 the Levi property. We shall write $L^0(\Phi)\!:=Z$ whenever $L^1(\Phi)$ is an order dense ideal
 in $Z$. Denote also by $L^\infty(\Phi)$ the order ideal in
 $Z$ generated by $\mathbb{1}$. Consider an order ideal
 $X\subset Z$ and we shall always assume
 that $L^\infty(\Phi)\subset X\subset L^1(\Phi)$. The associate
 space $X'$ is defined as the set of all $x'\in L^0(\Phi)$ for
 which $xx'\in L^1(\Phi)$ for all $x\in X$. Clearly, $X'$ is
 an order ideal in $Z$.

 If $(\Omega,\Sigma, \mu)$ is a probability space and  $\mathcal{X}_0$ is
 an~order closed vector sublattice of $L^\infty(\Omega,\Sigma, \mu)$
 containing $1_\Omega$, then there exists a $\sigma$-subalgebra $\Sigma_0$ of $\Sigma$
 such that $\mathcal{X}_0=L^\infty(\Omega,\Sigma_0,\mu_0)$,
 with $\mu_0=\mu|_{\mathcal{X}_0}$; cp.~\cite[Lemma
 2.2]{DHP}.

 Interpreting this fact and the properties of  conditional expectation in a Boolean value model yields the following result.

 \subsubsec{4.3.2.} \theorem{}Let $\Phi:L^1(\Phi)\to Y$ be a strictly
 positive Maharam operator with $Y=Y_\Phi$ and let $Z_0$ be an
 order closed sublattice in $L^0(\Phi)$.
 If $\mathbb{1}\in X_0\!:=L^1(\Phi)\cap Z_0$ and the restriction
 $\Phi_0\!:=\Phi|_{X_0}$ has the Maharam property then $X_0=L^1(\Phi_0)$
 and there exists an operator $\mathrm{E}(\cdot|Z_0)$ from
 $L^1(\Phi)$ onto $L^1(\Phi_0)$ such that
 \vspace*{2pt}

 \subsubsec{(1)}~$\mathrm{E}(\cdot|Z_0)$ is an~order
 continuous positive linear projection.
 \vspace*{2pt}

 \subsubsec{(2)}~$\mathrm{E}(\cdot|Z_0)$ commutes with all
 saturated projections, i.\,e. $\mathrm{E}(h(\pi)x|Z_0)=h(\pi)\mathrm{E}(x|Z_0)$
 for all $\pi\in \mathbb{P}_\Phi(X)$ and $x\in L^1(\Phi)$.
 \vspace*{2pt}


 \subsubsec{(3)}~$\Phi(xy)=\Phi(y\mathrm{E}(x|Z_0))$ for all $x\in L^1(\Phi)$
 and $y\in L^\infty(\Phi_0)$.
 \vspace*{2pt}

 \subsubsec{(4)}~$\Phi_0(|\mathrm{E}(x|Z_0)|)\leq\Phi(|x|)$
 for all $x\in L^1(\Phi)$.
 \vspace*{2pt}

 \subsubsec{(5)}~$\mathrm{E}(\cdot|Z_0)$ satisfies the averaging
 identity, i.\,e.
 $\mathrm{E}(v\mathrm{E}(x|Z_0)|Z_0)=\mathrm{E}(v|Z_0)\mathrm{E}(x|Z_0)$
 for all $x\in L^1(\Phi)$ and $v\in L^\infty(\Phi)$.
 \Endproc

 \subsubsec{4.3.3.}~We will call the operator $\mathrm{E}(\cdot|Z_0)$
 defined by Theorem 4.10.4 \textit{the conditional expectation
 operator with respect to\/} $Z_0$.
 Take $w\in X'$ and observe that $\mathrm{E}(wx|Z_0)\in
 L^1(\Phi_0)$ is well defined for all $x\in X$. If, moreover,
 $\mathrm{E}(wx|Z_0)\in X$ for every $x\in X$ then we can
 define a linear operator $T:X\to X$ by putting $Tx=\mathrm{E}(wx|Z_0)$
 $(x\in X)$. Clearly, $T$ is order bounded and order continuous.
 Furthermore, for all $x\in X_+$ we have
 $$
 T^+x=\mathcal{E}(w^+x|Z_0),\quad
 T^-x=\mathcal{E}(w^-x|Z_0),\quad
 |T|x=\mathcal{E}(|w|x|Z_0).
 $$
 In particular, $T$ is positive if and only if so is $w$. Putting
 $x\!:=wx$ and $y\!:=\mathbb{1}$ in 4.10.3\,(3), we get
 $\Phi(wx)=\Phi(wx\mathbb{1})=\Phi(\mathcal{E}(wx|Z_0))=\Phi(Tx)$
 for all $x\in X$. Now, $x$ can be chosen to be a component of
 $\mathbb{1}$ with $wx=w^+$ or $wx=w^-$, so that $T=0$ implies
 $\Phi(w^+)=0$ and $\Phi(w^-)=0$, since $\Phi$ is strictly positive.
 Thus $w\in X'$ is uniquely determined by $T$.

 Say that $T$ satisfies the \textit{averaging identity}, if
 $T(y\cdot Tx)=Ty\cdot Tx$ for all $x\in X$ and $y\in L^\infty(\Phi)$. Now we present two well-known results. By $\mathcal{E}(\cdot|\Sigma_0)$ we denote the conditional expectation operator with respect to $\sigma$-algebra $\Sigma_0$.

 \subsubsec{4.3.4.}~\proclaim{Theorem.}~Let $(\Omega,\Sigma,\mu)$ be probability
 space and $\mathcal{X}$ be an order ideal in $L^1(\Omega,\Sigma,\mu)$
 containing $L^\infty(\Omega,\Sigma,\mu)$. For a linear
 operator~$\mathcal{T}$ on~$\mathcal{X}$ the following are equivalent:

 \subsubsec{(1)}~$\mathcal{T}$ is order continuous, satisfies the averaging
 identity, and leaves invariant the subspace $L^\infty(\Omega,\Sigma,\mu)$.

 \subsubsec{(2)}~There exist $w\in\mathcal{X}'$ and a sub-$\sigma$-algebra
 $\Sigma_0$ of $\Sigma$ such that $\mathcal{T}x=\mathcal{E}(wx|\Sigma_0)$
 for all $x\in\mathcal{X}$.
 \Endproc

 The following two results can be proved by interpreting Theorems 4.3.4 and 4.3.5 in a Boolean valued model.

 \subsubsec{4.3.5.}~\proclaim{Theorem.}~For a subspace $\mathcal{X}$ of
 $L^1(\Omega,\Sigma,\mu)$ the following are equivalent:

 \subsubsec{(1)}~$\mathcal{X}$ is the range of a positive contractive
 projection.

 \subsubsec{(2)}~$\mathcal{X}$ is a closed vector sublattice
 of $L^1(\Omega,\Sigma,\mu)$.

 \subsubsec{(3)}~There exists a lattice isometry from some
 $L^1(\Omega',\Sigma',\mu')$ space onto $\mathcal{X}$.
 \Endproc

 \subsubsec{4.3.6.}~\proclaim{Theorem.}~Let $\Phi:L^1(\Phi)\to Y$ be a strictly
 positive Maharam operator and $X$ be an order dense ideal
 in $L^1(\Phi)$ containing $L^\infty(\Phi)$. For a linear
 operator $T$ on $X$ the following are equivalent:

 \subsubsec{(1)}~$T$ is order continuous, satisfies the averaging
 identity, leaves invariant the subspace $L^\infty(\Phi)$, and
 commutes with all $\Phi$-saturated projections.


 \subsubsec{(2)}~There exist $w\in X'$ and an order closed
 sublattice $Z_0$  in $L^0(\Phi)$ containing a unit element
 $\mathbb{1}$ of $L^1(\Phi)$ such that the restriction of
 $\Phi$ onto $L^1(\Phi)\cap Z_0$ has the Maharam property and
 $Tx=\mathcal{E}(wx|Z_0)$ for all $x\in X$. \Endproc

 \subsubsec{4.3.7.}~\proclaim{Theorem.}~For each subspace $X_0$ of $L^1(\Phi)$ the
 following statements are equivalent:

 \subsubsec{(1)}~$X$ is the range of a positive $\Phi$-contractive
 projection.

 \subsubsec{(2)}~$X$ is a closed vector sublattice of $L^1(\Phi)$
 invariant under all $\Phi$-saturated projections.

 \subsubsec{(3)}~There exists a Maharam operator $\Psi:L^1(\Psi)\to Y$
 and a lattice isomorphism $h$ from $L^1(\Psi)$ onto $X$
 such that $\Phi(|Tx|)=\Psi(|x|)$ for all $x\in L^1(\Psi)$.
 \Endproc

 \subsubsec{4.3.8.}~Theorems 4.3.4 and 4.3.6 can be found in \cite[Proposition~3.1]{DHP} and \cite[Lemma~1]{Dou1},
 respectively. Theorems 4.3.6 and 4.3.7 are published for the first time.

 \subsection{4.4. Maharam Extension}

 Thus, the general properties of Maharam operators can be deduced from the corresponding facts about functionals with the help of Theorem~4.1.8. Nevertheless, these methods  may be also useful in studying arbitrary regular operators.

 \subsubsec{4.4.1.}\proclaim{}Suppose that $X$ is a vector lattice over a dense subfield
 $\mathbb{F}\subset\mathbb{R}$ and $\varphi:X\to\mathbb{R}$ is a strictly positive
 $\mathbb{F}$-linear functional. There exist a Dedekind complete vector lattice
 $X^\varphi$ containing $X$ and a strictly positive order continuous linear functional
 $\bar{\varphi}:X^\varphi\to\mathbb{R}$ with the Levi property extending $\varphi$
 such that for every $x\in X^\varphi$
 there is a sequence $(x_n)$ in $X$ with $\lim_{n\to\infty}\bar{\varphi}(|x-x_n|)=0$.
 \Endproc

 \beginproof~Denote $d(x,y)\!:=\varphi(|x-y|)$ and note that $(X,d)$ is a metric space.
 Let $X^\varphi$ the completion of the metric space $(X,d)$ and $\bar{\varphi}$ is an extension
 of $\varphi$ to $X^\varphi$ by continuity. It is not difficult to ensure that $X^\varphi$ is
 a Banach lattice with an additive
 norm $\|\cdot\|^\varphi\!:=\bar{\varphi}(|\cdot|))$ containing $X$ as a norm dense
  $\mathbb{F}$-linear sublattice. Thus,  $\bar{\varphi}$ is a~strictly positive order continuous
  linear functional on $X^\varphi$ with the Levi Property.~\endproof

 \subsubsec{4.4.2.}~Denote $L^1(\varphi)\!:=X^\varphi$ and let $\bar{X}$ stands for the order
 ideal in $L^1(\varphi)$ generated by $X$. Then $(L^1(\varphi),\|\cdot\|^\varphi)$ is an
 $AL$-space and $\bar{X}$ is a Dedekind complete vector lattice. Moreover, $X$ is norm dense
 in $L^1(\varphi)$ and hence in $\bar{X}$.

 Given a nonempty subset $U$ of a lattice $L$,
 we denote by $U^{\scriptscriptstyle\uparrow}$ (resp. $U^{\scriptscriptstyle\downarrow}$)
 the set of elements
 $x\in L$ representable in the form $x=\sup(A)$ (resp. $x=\inf(A)$), where $A$ is an upward
 (resp. downward) directed subset of $U$. Moreover, we set  $U^{\scriptscriptstyle\uparrow\downarrow}\!:=
 (U^{\scriptscriptstyle\uparrow})^{\scriptscriptstyle\downarrow}$ etc. If in the above
 definition $A$ is countable, then we write $U^{\scriptscriptstyle\upharpoonleft}$,
 $U^{\scriptscriptstyle\downharpoonleft}$, and  $U^{\scriptscriptstyle\upharpoonleft\downharpoonleft}$
 instead of $U^{\scriptscriptstyle\uparrow}$, $U^{\scriptscriptstyle\downarrow}$, and
 $U^{\scriptscriptstyle\uparrow\downarrow}$. Recall that for the Dedekind completion $X^\delta$ we have $X^\delta=X^{\scriptscriptstyle\uparrow}=X^{\scriptscriptstyle\downarrow}$.

 \subsubsec{4.4.3.}\proclaim{}The identities
 $\bar{X}=X^{\scriptscriptstyle\downharpoonleft\upharpoonleft}
 =X^{\scriptscriptstyle\upharpoonleft\downharpoonleft}$ and
 $L^1(\varphi)=X^{\scriptscriptstyle\downharpoonleft\upharpoonleft}
 =X^{\scriptscriptstyle\upharpoonleft\downharpoonleft}$ hold with
 both
 $(\cdot)^{\scriptscriptstyle\downharpoonleft\upharpoonleft}$ and
 $(\cdot)^{\scriptscriptstyle\upharpoonleft\downharpoonleft}$
 taken in $\bar{X}$ and $L^1(\varphi)$, respectively.
 \Endproc

 Translating 4.4.1 and 4.4.2 by means of Boolean valued ``scalarization'' leads to the following result.

 \subsubsec{4.4.4.}~\proclaim{Theorem.}~Let $X$ and $Y$ be vector lattices with $Y$ Dedekind
 complete and $T$ a strictly positive linear operator from $X$ to $Y$. There exist a~Dedekind
 complete vector lattice $\bar{X}$ and a strictly positive Maharam operator $\bar{T}:\bar{X}\to Y$
 satisfying the conditions:

 \subsubsec{(1)}~There exist a lattice homomorphism $\iota:X
 \to\bar X $ and  an $f$-algebra homomorphism $\theta:\mathcal{Z}(Y)\to\mathcal{Z}(\bar X)$ such that
 \begin{equation}
 \alpha Tx=\bar{T}(\theta(\alpha)\iota(x))
 \quad(x\in X,\,\alpha\in\mathcal{Z}(Y)).
 \end{equation}

 \subsubsec{(2)}~$\iota(X)$ is a~majorizing sublattice in
 $\bar{X}$ and  $\theta(\mathcal{Z}(Y)) $ is an order closed sublattice and subring of $\mathcal{Z}(\bar{X})$.


 \subsubsec{(3)}~The representation
 $\bar{X}=(X\odot\mathcal{Z}(Y))^{\downarrow\uparrow}$ holds,
 where $X\odot\mathcal{Z}(Y)$ is a subspace of
 $\bar{X}$ consisting of all finite sums
 $\sum_{k=1}^n\theta(\alpha_k)\iota(x_k)$ with
 $x_1,\dots,x_n\in X$ and $\alpha_1,\dots,\alpha_n
 \in\mathcal{Z}(Y)$.
 \Endproc

 \subsubsec{4.4.5.}~The pair $(\bar{X}, \bar{T})$ (or $\bar{T}$ for short) is called a \textit{Maharam extension\/}
 of $T$ if it satisfies 4.4.4\,(1--3). The pair $(\bar{X},\iota)$
 is also called a \textit{Maharam extension space\/} for $T$. Two Maharam extensions $T_1$ and
 $T_2$ of $T$ with the respective Maharam extension spaces $(X_1,\iota_1)$ and $(X_2, \iota_2)$
 are said to be \textit{isomorphic\/} if there exists a lattice isomorphism $h$ of $X_1$ onto $X_2$
 such that $T_1=T_2\circ h$ and $\iota_2=h\circ\iota_1$. It is not difficult to ensure that
 a Maharam extension is unique up to isomorphism.

 \subsubsec{4.4.6.}~Let $X$ and $Y$ be vector lattices with $Y$ Dedekind complete,
 $T:X\to Y$ a strictly positive operator and $(\bar{X},\bar{T})$ a Maharam extension of $T$.
 Consider a universal completion $\bar{X}^\mathrm{u}$ of $\bar{X}$ with a fixed $f$-algebra structure.
 Let $L_1 (\Phi)$ be the greatest order dense ideal in $\bar{X}^\mathrm{u}$, onto which $\bar{T}$ can
 be extended by order continuity. In more details,
 $$
 \gathered
 L^1(T)\!:=\{x\in\bar{X}^\mathrm{u}:\ \bar{T}([0,\,|x|]\cap\bar{X})\text{~is order bounded in~} Y\},
 \\
 \hat{T}x\!:=\sup\{\bar{T}u:\ u\in\bar{X},\,0\leq u\leq x\} \quad(x\in L^1(T)_+),
 \\
 \hat{T}x=\hat{T}x^+-\hat{T}x^-\quad(x\in L^1(T)).
 \endgathered
 $$
 Define an $Y$-valued norm $[\!]\cdot[\!]$ on $L^1(T)$ by $[\!]u[\!]\!:=\hat{T}(|u|)$.
 In terminology of lattice normed spaces $(L^1(T),[\!]\cdot[\!])$ is a \textit{Banach--Kantorovich lattice},
 see \cite[Chapter 2]{DOP}. In~particular, $[\!]au[\!]=|a|[\!]u[\!]$ $(a\in\mathcal{Z}(Y),\,u\in L^1(T)$.

 \subsubsec{4.4.7.}~\theorem{}For every operator
 $S\in\{T\}^{\perp\perp}$ there is a~unique element
 $z=z_T\in \bar{X}^\mathrm{u}$ satisfying
 $$
 Sx =\hat{T}(z\cdot \imath(x))\quad (x \in X).
 $$
 The correspondence $T\mapsto z_T$ establishes a~lattice
 isomorphism between the band $\{T\}^{\perp\perp}$ and
 the~order dense ideal in $\bar{X}^\mathrm{u}$ defined by
 $$
 \{z\in\bar{X}^\mathrm{u}:\, z\cdot \imath(X)\subset L_1(T)\}.
 $$
 \Endproc

 \beginproof~This result is a~variant of the Radon--Nikod\'ym Theorem for positive operators and may be obtained as a combination of Theorems 4.2.10 and 4.4.4 or proved by means of Boolean valued ``scalarization''.~\endproof

 \subsubsec{4.4.8.}~Maharam extension stems from the corresponding extension result given by D.~Maharam for $F$-integrals \cite{Mah3, Mah4, Mah5}. For operators in Dedekind complete vector lattices this construction was performed in \cite{AKK, AKK1} by three different ways. One of them, based upon the imbedding $x\mapsto\bar x$ of a vector lattice $X$ into $L^\sim\big((L^\sim(X,Y),Y\big)$ defined as   $\bar x(T)\!:=Tx$ $(T\in L^\sim(X,Y))$, was independently discovered in \cite{LuP}. The main difference is is that in \cite{LuP} the Maharam extension was constructed for an arbitrary collection of order bounded operators. For some further properties of Maharam extension see in~\cite{DOP} and \cite{LuP}.


 \section{References}
 \def\citeitem#1#2{\bibitem{#2}}

 \begin{enumerate}
 {\itemsep=0pt\parskip=0pt
 \leftskip=-10pt
 \normalsize

 \citeitem{AbramovichVK1979}{AVK} \textsl{Abramovich Yu.~A., Veksler A.~I. and Koldunov A.~V.}
 On~dis\-joint\-ness preserv\-ing operators~/\!/ Dokl. Akad. Nauk SSSR.---1979.---Vol.~289, No.~5.---pp.~1033--1036.

 \citeitem{AbramovichVK1981}{AVK1} \textsl{Abramovich Yu.~A., Veksler A.~I., and Koldunov A.~V.} 
 Disjointness preserving
 operators, their continuity, and multiplicative representation~/\!/ Linear Operators and Their Appl.---Leningrad:
 Leningrad Ped. Inst., 1981.---[in Russian].

 \citeitem{AbramovichKi1999}{AK1} \textsl{Abramovich Y.~A.and Kitover A.~K.} $d$-Independence and $d$-bases in Vector lattices~/\!/
 Rev. Romaine de Math. Pures et Appl.---1999.---Vol.~44.---pp.~667--682.

 \citeitem{AbramovichKi2000}{AK3} \textsl{Abramovich Y. A. and Kitover A.~K.} Inverses of disjointness preserving operators~/\!/
 Memoirs Amer. Math. Soc.---2000.---Vol.~679.

 \citeitem{AbramovichKi2003}{AK2} \textsl{Abramovich Y.~A.and Kitover~A.~K.} $d$-Independence and $d$-bases~/\!/
 Positivi\-ty.---2003.---Vol.~7, No.~1.---pp.~95--97.

 \citeitem{AczelDhombres}{AD}      \textsl{Acz\'el J. and Dhombres J.} Functional Equations in Several
 Variables.---Cambrid\-ge etc.: Cambridge Univ. Press, 1989.

 \citeitem{AkilovKut1978}{AK}      \textsl{Akilov G.~P. and Kutateladze S.~S.} Ordered Vector  Spaces.---Novosibirsk: Nauka, 1978.

 \citeitem{Akilov1988}{AKK}     \textsl{Akilov G.~P., Kolesnikov E.~V., and Kusraev A.~G.} The Lebesgue
 extension of a positive operator~/\!/ Dokl. Akad. Nauk SSSR.---1988.---Vol.~298,
 No.~3.---pp.~521--524.---[in Russian].

 \citeitem{AkilovіKol1988}{AKK1}    \textsl{Akilov G.~P., Kolesnikov E.~V., and Kusraev A.~G.}
 On order continuous extension of a positive operator~/\!/ Sibirsk. Mat. Zh.---1988.---Vol.~29, No.~5.---pp.~24--35.---[in Russian].

 \citeitem{AliprantisBur1985}{AB}   \textsl{Aliprantis C.~D. and Burkinshaw~O.} Positive Operators.---N.\,Y.: Acad. Press, 1985.---367~p.

 \citeitem{Bell~J.~L.}{Bell}    \textsl{Bell~J.~L.} Boolean-Valued Models and Independence
 Proofs in Set Theory.---N.\,Y.: Clarendon Press, 1985.---xx+165~p.

 \citeitem{BernauHuijsmans}{BHP} \textsl{Bernau~S.~J., Huijsmans C.~B., and de Pagter B.} Sums of lattice
 homomorphisms~/\!/ Proc. Math. Soc.---1992.---Vol.~115, No.~1.

 \citeitem{Boulabiar}{Boul} \textsl{Boulabiar K.} Recent trends on order bounded disjointness preserving
 operators~/\!/ Irish Math. Soc. Bulletin.---2008.---Vol.~62.---pp.~43--69.

 \citeitem{BoulabiarSS}{BBS} \textsl{Boulabiar K., Buskes G., and Sirotkin G.} Algebraic order bounded
 disjointness preserving operators and strongly diagonal operators~/\!/
 Integral Equations and Operator Theory.---2006.---Vol.~54.---pp.~9--31.

 \citeitem{Bourgain1986}{Bourg} \textsl{Bourgain J.} Real isomorphic complex Banach spaces need not
 be complex isomorphic~/\!/ Proc. Amer. Math. Soc.---1986.---Vol.~96.---pp.~221--226.

 \citeitem{Cohen1966}{Coh}    \textsl{Cohen~P.~J.} Set Theory and the
 Continuum Hypothesis.---N.\,Y.--Amster\-dam: W.~A.~Benjamin Inc., 1966.

 \citeitem{Cooper}{Coop}  \textsl{Cooper J.~L.~B.} Coordinated linear spaces~/\!/
 Proc. London Math. Soc.---1953.---Vol.~3.---pp.~305--327.

 \citeitem{Dieudonne1952}{Die}  \textsl{Dieudonn\'e~J.} Complex structures on real Banach spaces~/\!/
 Proc. Amer. Math. Soc.---1952.---Vol.~3.---pp.~162--164.

 \citeitem{Dodds}{DHP} \textsl{Dodds~P.~G.,  Huijsmans~C.~B., and de Pagter~B.} Characterizations
 of conditional expectation-type operators~/\!/ Pacific J.  Math.---1990.---Vol.~141, No.~1.

 \citeitem{Douglas1982}{Dou1}  \textsl{Douglas R.} Contractive projections on $L_1$ space~/\!/
 Pacific J. Math. 1965.---Vol.~15, No.~2.---pp.~443--462.

 \citeitem{DuhouxMeyer1982}{DuM}  \textsl{Duhoux~M. and Meyer~M.} A new proof of the lattice
 structure of orthomorphisms~/\!/ J. London Math. Soc.---1982.---Vol.~25, No.~2.---pp.~375--378.

 \citeitem{Gordon1977}{Gor1}   \textsl{Gordon~E.~I.}
 Real numbers in Boolean-valued models of set theory and $K$-spa\-ces~/\!/
 Dokl. Akad. Nauk SSSR.---1977.---Vol.~237, No.~4.---pp.~773--775.

 \citeitem{Gordon1981}{Gor2}  \textsl{Gordon~E.~I.}
 $K$-spaces in Boolean-valued models of set theory~/\!/
 Dokl. Akad. Nauk SSSR.---1981.---Vol.~258, No.~4.---pp.~777--780.

 \citeitem{Gordon1982}{Gor3}  \textsl{Gordon~E.~I.}
 To the theorems of identity preservation  in~$K$-spaces~/\!/ Sibirsk. Mat. Zh.---1982.---Vol.~23, No.~5.---pp.~55--65.

 \citeitem{Gowers1993}{Gow}  \textsl{Gowers W.~T.} A solution to Banach's hyperplane problem~/\!/ Bull.
 London Math. Soc.---1994.---Vol.~26.---pp.~523--530.

 \citeitem{GowersMaur1994}{GM1}  \textsl{Gowers~W. T. and Maurey~B.} The unconditional basic
 sequence problem~/\!/ J. Amer. Math. Soc.---1993.---Vol.~6.---pp.~851--874.

 \citeitem{Gowers1994}{GM2}  \textsl{Gowers~W. T. and Maurey~B.} Banach spaces with small
 spaces of operators~/\!/  Math. Ann.---1997.---Vol.~307.---pp.~543--568.

 \citeitem{Grothendieck}{Gro}  \textsl{Grothendieck~A.}
 Une caract\'erisation vectorielle-m\'etrique des espae $L^1$~/\!/
 Canad. J. Math.---1955.---Vol.~4.---pp.~552--561.

 \citeitem{Gutman1995}{Gut6}  \textsl{Gutman A.~E.} Locally one-dimensional $K$-spaces and $\sigma$-distributi\-ve
 Boole\-an algebras~/\!/ Siberian Adv. Math.---1995.---Vol.~5, No.~2.---pp.~99--121.

 \citeitem{Gutman1996}{Gut1}  \textsl{Gutman~A.~E.} Disjointness preserving operators~/\!/ Vector Lattices and Integral Operators
 (Ed. S.~S.~Kutateladze).---Dordrecht etc.: Kluwer,  1996.---pp.~361--454.

 \citeitem{GutmanKK}{GKK}  \textsl{Gutman~A.~E., Kusraev A.~G., and Kutateladze S.~S.}
 The Wickstead Problem. ---Vladikavkaz, 2007.---44~p.---(Preprint~/
 IAMI VSC RAS; No.~3).

 \citeitem{HuijsmansWic1992}{HW}  \textsl{Huijsmans~C.~B. and Wickstead A.~W.} The inverse
 of band preserving and disjointness preserving operators~/\!/
 Indag. Math.---1992.---Vol.~3, No.~2.---pp.~179--183.

 \citeitem{Jech97}{Jech} \textsl{Jech~T.~J.} Set Theory.---Berlin: Springer-Verlag, 1997.---634~p.


 \citeitem{KantorovichA1984}{KA}   \textsl{Kantorovich~L. V. and Akilov~G.~P.} Functional Analysis.---Moscow:
 Nauka, 1984.---[in Russian].

 \citeitem{KantorovichVP1950}{KVP} \textsl{Kantorovich~L. V., Vulikh~B. Z., andPinsker~A. G.} Functional
 Analysis in~Semi\-ord\-er\-ed Spaces.---M.--L.: Gostekhizdat, 1950.---[in Russian].

 \citeitem{Krasnosel1985}{KLS}     \textsl{Krasnosel'ski\u\i{}~M.~A.,  Lifshits~E.~A., and Sobolev~A.~V.}
 Positive Linear Systems. The Method of~Positive  Operators---Moscow:
 Nauka, 1985.---[in~Russian].

 \citeitem{Kuczma}{Kuc}  \textsl{Kuczma M.} An introduction to the theory of functional equations and
 in\-equ\-a-li\-ties.---Basel--Boston--Berlin: Birkh\"auser, 2009.

 \citeitem{Kusraev1982a}{K12}     \textsl{Kusraev A. G.} General desintegration formulas~/\!/
 Dokl. Akad. Nauk SSSR.---1982.---Vol.~265, No.~6.---pp.~1312--1316.

 \bibitem{Kus1}  \textsl{Kusraev A.~G.} Order continuous
 functionals in Boolean-valued models
 of set theory~/\!/ Siberian Math. J.---1984.---Vol.~25, No.~1.---pp.~57--65.

 \citeitem{Kusraev1986}{K9}  \textsl{Kusraev~A.~G.} Numeric systems in Boolean-valued models of set
 theory~/\!/ Proceed. of the VIII All-Union Conf. on Math. Logic (Moscow).---Moscow, 1986.---P.~99.

 \citeitem{Kusraev2000}{DOP}   \textsl{Kusraev~A.~G.} Dominated Operators.---Dordrecht: Kluwer, 2000.

 \citeitem{Kusraev2003}{DOR}   \textsl{Kusraev~A.~G.} Dominated Operators.---Moscow, Nauka, 2003.---619~p.

 \citeitem{Kusraev2004}{Kus11}   \textsl{Kusraev A.~G.} On band preserving operators~/\!/  Vladikavkaz Math.
 J.---2004.---Vol.~6, No.~3.---pp.~47--58.---[in Russian].

 \citeitem{Kusraev2006}{Kus2}  \textsl{Kusraev~A.~G.}
 Automorphisms and derivations in extended complex $f$-algebras~/\!/ Siberian Math. J.---2006.---Vol.~47, No.~1.---pp.~97--107.

 \citeitem{Kusraev2006}{Kus4}  \textsl{Kusraev~A.~G.} Analysis, algebra, and logics in operator theory~/\!/
 Complex Analysis, Operator Theory, and Mathematical Modeling
 (Eds. Korobe\u{\i}nik~Yu.~F., Kusraev A.~G.).---Vladikavkaz: VSC RAS, 2006.---pp.~171--204.---[in Russian].

 \citeitem{KusKut1992}{KK3}   \textsl{Kusraev A. G. and Kutateladze S. S.} Subdifferentials:
 Theory and Appli\-cati\-ons.---Novosibirsk: Nauka, 1992.---270~p.;
 Dordrecht: Kluwer Academic Pub\-li\-shers, 1995.---398~p.

 \citeitem{KusraevKut1990}{KK}    \textsl{Kusraev~A. G. and Kutateladze~S. S.} Nonstandard Methods of
 Analysis.---Novosibirsk: Nauka,  1990; Dordrecht: Kluwer Academic Publisher, 1994.

 \citeitem{KusraevKut1999}{BVA} \textsl{Kusraev~A. G. and Kutateladze~S. S.} Boolean-Valued
 Analysis.---Novosibirsk: Nauka, 1999; Dordrecht: Kluwer Academic Publisher, 1999.

 \citeitem{KusraevKut2005}{IBA}    \textsl{Kusraev~A.~G. and Kutateladze~S.~S.} Introduction
 to Boolean-Valued Analy\-sis.---M.: Nauka, 2005.---[in Russian].

 \citeitem{KusraevKut2014}{KK_DPO}
 \textsl{Kusraev,~A.~G. and Kutateladze,~S.~S.,}
On order bounded disjointness preserving operators /\!/ Sib. Math. J. 2014. V.~55, N~5. 915--928.

 \citeitem{KusraevKur_selected}{KK_T} \textsl{Kusraev~A. G. and Kutateladze~S. S.}
 Boolean Valued Analy\-sis: Selected Topics. Vladikavkaz:
 SMI VSC RAS, 2014. (Trends in Science: The South of Russia. A~Mathematical Monograph. Issue~6.)

 \citeitem{KusraevTab}{KT1}  \textsl{Kusraev~A.~G. and  Tabuev~S.~N.} On
 multiplicative representation of disjointness preserving bilinear operators~/\!/
 Sib. Math. J. 2008. V.~49, N~2. P. 357--366.

 \citeitem{Kusraeva}{KZ} \textsl{Kusraeva~Z.~A.} Band preserving algebraic operators~/\!/
 Vladikavkaz Math. J.---2013.---Vol.~15, No.~3.---pp.~54--57.

 \citeitem{Kusraeva}{KZ1} \textsl{Kusraeva Z.~A.} Involutions and
 complex structures on real vector lattices.---To appear.

 \citeitem{Kutateladze1979}{Kut1}  \textsl{Kutateladze~S.~S.} Choquet boundaries in $K$-spaces~/\!/ Uspekhi Mat. Nauk.---1975.---Vol.~30, No.~4.---pp.~107--146.

 \citeitem{Kutateladze1979}{Kut2}  \textsl{Kutateladze~S.~S.} Convex operators~/\!/ Uspekhi Mat. Nauk.---1979.---Vol.~34, No.~1.---pp.~167--196.

 \citeitem{Kutateladze2005}{KutH}    \textsl{Kutateladze~S.~S.} On differences of lattice
 homomorphisms /\!/ Siberian Math. J.---2005.---Vol.~46, No.~2.---pp.~393--396.

 \citeitem{Kutateladze2005}{KutG}    \textsl{Kutateladze~S.~S.} On Grothendieck subspaces~/\!/
 Siberian Math. J.---2005.---Vol.~46, No.~3.---pp.~620--624.

 \citeitem{Lacey}{Lac}             \textsl{Lacey~H.~E.} The Isometric Theory of Classical Banach
 Spaces.---Berlin etc.: Springer-Verlag,  1974.---247~p.

 \citeitem{LindenstraussWul}{LW} {\sl Lindenstrauss J. and Wulbert~D.~E.} On the classification of
 the Banach spaces whose duals are $L^1$-spaces~/\!/ J. Funct. Anal.---1969.---Vol.~4.---pp.~332--349.

 \citeitem{LuxemburgPagt}{LuP}     \textsl{Luxemburg~W.~A.~J. and de Pagter~B.} Maharam extension
 of positive operators and $f$-algebras~/\!/ Positivity.---2002.---Vol.~6, No.~2.---pp.~147--190.

 \citeitem{LuxemburgSchep}{LS}    \textsl{Luxemburg~W.~A.~J. and Schep~A.}
 A Radon-Nikod\'ym type theorem for positive operators and a dual~/\!/ Indag. Math.---1978.---Vol.~40.---pp.~357--375.

 \citeitem{LuxZa}{LuZ} \textsl{Luxemburg~W.~A.~J. and Zaanen~A.~C.} Riesz Spaces. Vol.~1.---Amster\-dam--London:
 North-Holland, 1971.---514~p.

 \citeitem{Mah53}{Mah3} \textsl{Maharam~D.} The representation of abstract integrals~/\!/
 Trans. Amer. Math. Soc.---1953.---Vol.~75, No.~1.---pp.~154--184.

 \citeitem{Mah55}{Mah4} \textsl{Maharam~D.}  On kernel representation of linear operators~/\!/
 Trans. Amer. Math. Soc.---1955.---Vol.~79, No.~1.---pp.~229--255.

 \citeitem{Mah84}{Mah5} \textsl{Maharam~D.} On positive operators~/\!/ Contemp. Math.---1984.---Vol.~26.---pp.~263--277.

 \citeitem{McPolinWick1985}{MW}    \textsl{McPolin~P.~T.~N. and Wickstead~A.~W.} The order boundedness
 of band preserving operators on uniformly complete vector lattices~/\!/
 Math. Proc. Cambridge Philos. Soc.---1985.---Vol.~97, No.~3.---pp.~481--487.

 \citeitem{Mey76}{Mey1} \textsl{Meyer~M.} Le stabilisateur d'un espace vectoriel r\'eticul\'e~/\!/
 C. R. Acad. Sci. Ser.~A.---1976.---Vol.~283.---pp.~249--250. 

 \citeitem{MeyNieb}{MN} \textsl{Meyer-Nieberg~P.} Banach Lattices.---Berlin etc.: Springer-Verlag, 1991.---xv+395~p.


 \citeitem{Schaefer1974}{Sch}  \textsl{Schaefer~H.~H.}
 Banach Lattices and Positive Operators.---Berlin etc.: Sprin\-ger-Verlag,  1974.---376~p.

 \citeitem{Schwarz1984}{Schw} \textsl{Schwarz~H.-V.} Banach Lattices and Operators.---Leipzig: Teubner, 1984.---208~p.

 \citeitem{Sem71}{Sem} \textsl{Semadeni~Zb.} Banach Spaces of Continuous Functions.---War\-szawa:
 Polish Scientific Publ., 1971.---584~p.

 \citeitem{Sikorski1964}{Sik}  \textsl{Sikorski R.} Boolean
 Algebras.---Berlin etc.: Springer-Verlag, 1964.

 \citeitem{SolovayTen}{SoT} \textsl{Solovay~R. and Tennenbaum~S.} Iterated Cohen extensions and
 Souslin's problem~/\!/ Ann. Math.---1972.---Vol.~94, No.~2.---pp.~201--245.

 \citeitem{Szarek}{Sza}    \textsl{Szarek~S.} A superreflexive Banach space which does not admit
 complex structure~/\!/ Proc. Amer. Math. Soc.---1986.---Vol.~97.---pp.~437--444.

 \citeitem{Tak1978}{Tak}   \textsl{Takeuti~G.} Two Applications of Logic to Mathematics.---Tokyo
 and Princeton: Iwanami and Princeton Univ. Press,
 1978.---137~p.

 \citeitem{Tak79}{Take} \textsl{Takeuti~G.} A transfer principle in harmonic analysis~/\!/
 J. Symbolic Logic.---1979.---Vol.~44, No.~3.---pp.~417--440.

 \citeitem{Takeuti1979}{Tak1} \textsl{Takeuti~G.} Boolean-valued analysis~/\!/ Applications
 of Sheaves (Proc. Res. Sympos. Appl. Sheaf Theory to Logic, Algebra and Anal., Univ. Durham,
 Durham, 1977).---Berlin etc.: Springer-Verlag, 1979. ---pp.~714--731.---(Lec\-ture Notes in Math. Vol.~753).

 \bibitem{TZ}   \textsl{Takeuti~G. and Zaring~W.~M.} Axiomatic set Theory.---N.\,Y.: Springer-Verlag, 1973.---238~p.

 \citeitem{Vladimirov1969}{Vl}    \textsl{Vladimirov~D.~A.} Boolean Algebras.---Moscow: Nauka,
 1969.---[in Russian].

 \citeitem{Vulikh1961}{Vul}       \textsl{Vulikh~B.~Z.} Introduction to the Theory of Partially Ordered
 Spaces.---Moscow: Fizmatgiz, 1961.

  \citeitem{VulikhLoz}{VL}           \textsl{Vulikh B.~Z. and  Lozanovskii G. Ya.} On the representation of completely
 linear and regular functionals in partially ordered spaces~/\!/ Mat. Sb.---1971.---Vol.~84\,(126),
 No.~3.---pp.~331--352.

 \citeitem{Wickstead1977}{Wic1}   \textsl{Wickstead~A.~W.} Representation and duality of multiplication
 operators on Archimedean Riesz spaces~/\!/ Compositio Math.---1977.---Vol.~35, No.~3.---pp.~225--238.

 \citeitem{Zaanen1983}{Z}  \textsl{Zaanen~A.~C.} Riesz Spaces. 2.---Amsterdam etc.:
 North-Holland, 1983.---720~p.

}
\end{enumerate}

\noindent
{\it Sem\"en S.~Kutateladze}
\par\smallskip
{\leftskip\parindent\small
\noindent
Sobolev Institute of Mathematics\\
4 Koptyug Avenue\\
Novosibirsk, 630090, RUSSIA\\
E-mail: sskut@math.nsc.ru
\par}

\noindent
{\it Anatoly G.~Kusraev}
\par\smallskip
{\leftskip\parindent\small
\noindent
Southern Mathematical Institute\\
 22 Markus Street\\
Vladikavkaz, 362027, RUSSIA\\
E-mail: kusraev@smath.ru
\par}

 \end{document}